%% file: HT.tex
\documentclass[12pt]{amsart}

\usepackage{latexsym, amssymb, amscd, amsfonts}

\usepackage{microtype}
\usepackage{pstricks}
\usepackage{multido}
\usepackage{pst-node}
\usepackage{pst-arrow}
\usepackage{pst-plot} 
\usepackage[colorlinks=true,urlcolor=black,citecolor=black,linkcolor=black]{hyperref}
\usepackage{graphicx}
\usepackage[all]{xypic}
\usepackage{ifthen}

\newboolean{bourbaki}
\setboolean{bourbaki}{true}

\newboolean{draft}
\setboolean{draft}{false}

\newboolean{shortnote}
\setboolean{shortnote}{false}

\ifthenelse{\boolean{shortnote}}{%
\newcommand{\point}{\vspace{3mm}\par\refstepcounter{subsection}\noindent{\bf \thesubsection.} }

\newcommand{\bpoint}[1]{\vspace{3mm}\par\refstepcounter{subsection}\noindent{\bf \thesubsection.} 
  {\bf #1.} }

\renewenvironment{equation}{\medskip\noindent\refstepcounter{subsection}\makebox[0pt][l]{({\bf\thesubsection})}\begin{minipage}{\textwidth}$$}{$$\end{minipage}\medskip\noindent}
} 
{
\newcommand{\point}{\vspace{3mm}\par\refstepcounter{subsection}\noindent{\bf \thesubsection.} }

\newcommand{\bpoint}[1]{\vspace{3mm}\par\refstepcounter{subsection}\noindent{\bf \thesubsection.} 
  {\bf #1.} }

\renewenvironment{equation}{\medskip\noindent\refstepcounter{subsubsection}\makebox[0pt][l]{({\bf\thesubsubsection})}\begin{minipage}{\textwidth}$$}{$$\end{minipage}\medskip\noindent}

}

\newcommand{\Fig}{\refstepcounter{subsubsection}{\sc Figure \thesubsubsection.}} 

\newcommand{\bpf}{\noindent {\em Proof.  }}
\newcommand{\epf}{\qed \vspace{+10pt}}

\catcode`\@=11 
\def\addtopunct#1{\expandafter\let\csname punct@\meaning#1\endcsname\let}
\addtopunct{.}    \addtopunct{,}    \addtopunct{?}
\addtopunct{!}    \addtopunct{;}    \addtopunct{:}

\let\seveendformula\]
\def\PunctAndEndFormula #1{#1\seveendformula}
\def\]{\futurelet\punctlet\checkpunct@i}
\def\checkpunct@i{\expandafter\ifx\csname punct@\meaning\punctlet\endcsname\let  
       \expandafter\PunctAndEndFormula 
       \else \expandafter\seveendformula\fi}
\catcode`\@=12 

\makeatletter
\newcommand*{\da@rightarrow}{\mathchar"0\hexnumber@\symAMSa 4B }
\newcommand*{\da@leftarrow}{\mathchar"0\hexnumber@\symAMSa 4C }
\newcommand*{\xdashrightarrow}[2][]{%
  \mathrel{%
    \mathpalette{\da@xarrow{#1}{#2}{}\da@rightarrow{\,}{}}{}%
   }%
}
\newcommand{\xdashleftarrow}[2][]{%
  \mathrel{%
    \mathpalette{\da@xarrow{#1}{#2}\da@leftarrow{}{}{\,}}{}%
  }%
}
\newcommand*{\da@xarrow}[7]{%
  \sbox0{$\ifx#7\scriptstyle\scriptscriptstyle\else\scriptstyle\fi#5#1#6\m@th$}%
  \sbox2{$\ifx#7\scriptstyle\scriptscriptstyle\else\scriptstyle\fi#5#2#6\m@th$}%
  \sbox4{$#7\dabar@\m@th$}%
  \dimen@=\wd0 %
  \ifdim\wd2 >\dimen@
    \dimen@=\wd2 %
  \fi
  \count@=2 %
  \def\da@bars{\dabar@\dabar@}%
  \@whiledim\count@\wd4<\dimen@\do{%
    \advance\count@\@ne
    \expandafter\def\expandafter\da@bars\expandafter{%
      \da@bars
      \dabar@ 
    }%
  }%
  \mathrel{#3}%
  \mathrel{%
    \mathop{\da@bars}\limits
    \ifx\\#1\\%
    \else
      _{\copy0}%
    \fi
    \ifx\\#2\\%
    \else
      ^{\copy2}%
    \fi
  }%
  \mathrel{#4}%
}
\makeatother

\newcommand{\xmapsto}{\mapstochar\relbar\joinrel\xrightarrow}

\newcommand{\xdownarrow}[1]{%
  {\left\downarrow\vbox to #1{}\right.\kern-\nulldelimiterspace}
}

\makeatletter
\newcommand*{\Relbarfill@}{\arrowfill@\Relbar\Relbar\Relbar}
\newcommand*{\xeq}[2][]{\ext@arrow 0055\Relbarfill@{#1}{#2}}
\makeatother

\newgray{lgray}{0.9}
\newgray{vlgray}{0.95}
\newgray{midgray}{0.65}
\newgray{darkgray}{0.3}

\newcommand{\Scale}{%
dup 4 -1 roll mul 3 1 roll mul
}

\newcommand{\SetTT}[6]{%
\def \TT {[ #1 #2 #3 #4 6 -1 roll dup #5 mul exch #6 mul ] transform }
}
\SetTT{-0.5}{-0.5}{1}{0}{0}{1}

\newcommand{\RotAngleSetTT}[2]{%
\SetTT{#1 cos neg}{#1 sin #2 sin mul neg}{#1 sin}{#1 cos #2 sin mul neg}{0}{#2 cos}
}

\newcommand{\TTrans}{%
4 -1 roll add exch 4 -1 roll add exch 4 2 roll add 3 1 roll 
}

\newcommand{\TScale}{%
 dup 5 -1 roll mul 4 1 roll \Scale 
}

\newcommand{\MakeUnit}{%
 3 copy dup mul 3 1 roll dup mul 3 1 roll dup mul add add sqrt 1 exch div \TScale 
}

\newcommand{\Hom}{\operatorname{Hom}}

\renewcommand{\part}{\operatorname{part}}

\newcommand{\st}{\,\,|\,\,}

\newcommand{\Fsh}{\mathcal{F}}

\newcommand{\Ish}{\mathcal{I}}
\newcommand{\Osh}{\mathcal{O}}
\newcommand{\Zsh}{\mathcal{Z}}
\renewcommand{\geq}{\geqslant}
\renewcommand{\leq}{\leqslant}

\renewcommand{\emptyset}{\varnothing} 
\renewcommand{\phi}{\varphi}

\newcommand{\Spec}{\operatorname{Spec}} 

\newcommand{\sg}{\mbox{\tiny sg}}   
\newcommand{\charx}{\mathfrak{X}}   
\newcommand{\expi}{\operatorname{expi}} 
\newcommand{\Flag}{F}   
\renewcommand{\Im}{\operatorname{Im}}   
\newcommand{\Bsub}{v}   

\makeatletter
\newcommand{\smallbullet}{} 
\DeclareRobustCommand\smallbullet{%
  \mathord{\mathpalette\smallbullet@{0.65}}
}
\newcommand{\smallbullet@}[2]{%
  \vcenter{\hbox{\scalebox{#2}{$\m@th#1\bullet$}}}%
}
\makeatother

\ifthenelse{\boolean{bourbaki}}
{
\renewcommand{\AA}{\mathbf{A}} 
\newcommand{\CC}{\mathbf{C}} 
\newcommand{\PP}{\mathbf{P}} 
\newcommand{\QQ}{\mathbf{Q}} 
\newcommand{\RR}{\mathbf{R}} 
\newcommand{\ZZ}{\mathbf{Z}} 

\DeclareSymbolFont{frenchmath}{OT1}{cmr}{m}{n} 

\DeclareMathSymbol{A}{\mathalpha}{frenchmath}{`A}  
\DeclareMathSymbol{B}{\mathalpha}{frenchmath}{`B}  
\DeclareMathSymbol{C}{\mathalpha}{frenchmath}{`C}
\DeclareMathSymbol{D}{\mathalpha}{frenchmath}{`D}
\DeclareMathSymbol{E}{\mathalpha}{frenchmath}{`E}
\DeclareMathSymbol{F}{\mathalpha}{frenchmath}{`F}
\DeclareMathSymbol{G}{\mathalpha}{frenchmath}{`G}
\DeclareMathSymbol{H}{\mathalpha}{frenchmath}{`H}
\DeclareMathSymbol{I}{\mathalpha}{frenchmath}{`I}
\DeclareMathSymbol{J}{\mathalpha}{frenchmath}{`J}
\DeclareMathSymbol{K}{\mathalpha}{frenchmath}{`K}
\DeclareMathSymbol{L}{\mathalpha}{frenchmath}{`L}
\DeclareMathSymbol{M}{\mathalpha}{frenchmath}{`M}
\DeclareMathSymbol{N}{\mathalpha}{frenchmath}{`N}
\DeclareMathSymbol{O}{\mathalpha}{frenchmath}{`O}
\DeclareMathSymbol{P}{\mathalpha}{frenchmath}{`P}
\DeclareMathSymbol{Q}{\mathalpha}{frenchmath}{`Q}
\DeclareMathSymbol{R}{\mathalpha}{frenchmath}{`R}
\DeclareMathSymbol{S}{\mathalpha}{frenchmath}{`S}
\DeclareMathSymbol{T}{\mathalpha}{frenchmath}{`T}
\DeclareMathSymbol{U}{\mathalpha}{frenchmath}{`U}
\DeclareMathSymbol{V}{\mathalpha}{frenchmath}{`V}
\DeclareMathSymbol{W}{\mathalpha}{frenchmath}{`W}
\DeclareMathSymbol{X}{\mathalpha}{frenchmath}{`X}
\DeclareMathSymbol{Y}{\mathalpha}{frenchmath}{`Y}
\DeclareMathSymbol{Z}{\mathalpha}{frenchmath}{`Z}
} 
{ 
\renewcommand{\AA}{\mathbb{A}} 
\newcommand{\CC}{\mathbb{C}} 
\newcommand{\PP}{\mathbb{P}} 
\newcommand{\QQ}{\mathbb{Q}} 
\newcommand{\RR}{\mathbb{R}} 
\newcommand{\ZZ}{\mathbb{Z}} 
}

\ifthenelse{\boolean{draft}}
{
\newcommand{\remind}[1]{{\bf[#1]}}
\newcommand{\lremind}[1]{{\bf[label:  #1]}}
\newcommand{\bremind}[1]{{\bf[label:  #1]}}

}
{
\newcommand{\remind}[1]{{}}
\newcommand{\lremind}[1]{{}}
\newcommand{\bremind}[1]{{}}

}

\DeclareMathAlphabet{\matheuler}{U}{zeur}{m}{n}
\DeclareMathAlphabet{\matholdstyle}{OT1}{pplj}{m}{n}
\newcommand{\euler}[1]{\ensuremath{\matheuler{#1}}}  
\newcommand{\Euler}{\euler} 

\newcounter{enumcount}
\newcommand{\PauseEnumerate}{\end{enumerate}\setcounter{enumcount}{\arabic{enumi}}}
\newcommand{\ResumeEnumerate}{\begin{enumerate}\setcounter{enumi}{\theenumcount}}

\newcommand{\AlphaList}{\renewcommand{\labelenumi}{{({\em\alph{enumi}})}}}

\newcommand{\RomanList}{\renewcommand{\labelenumi}{{({\em\roman{enumi}})}}}
\newcommand{\EulerList}{\renewcommand{\labelenumi}{{($\matheuler{\arabic{enumi}}$)}}}
\AlphaList  

\begin{document}
\pagestyle{plain} \title{\large Homeomorphism type of the non-negative part of a complete toric variety}
\author{Mike Roth}
\address{Dept.\ of Mathematics and Statistics, Queens University, Kingston,
Ontario, Canada} 
\email{mike.roth@queensu.ca} 
\thanks{Research partially supported by an NSERC grant}

\begin{abstract}
In this note we show that the nonnegative part of a proper complex toric variety has the homeomorphism
type of a sphere, and consequently that the nonnegative part has a natural structure of a cell complex.  
This extends previous results of Ehlers and Jurkiewicz.
The proof also provides a simplicial decomposition of the nonnegative part, and
a parameterization of each maximal simplex.
\end{abstract}

\maketitle

\section{Introduction}

\point
\label{sec:non-neg-part-def}
Let $X$ be an $n$-dimensional toric variety over $\CC$ with fan $\Delta_{X} \subseteq N_{\RR}$, 
where $M$ and $N$ denote as usual the lattice of torus characters and one-parameter subgroups respectively.  

Since the transition data between the torus-stable open affine
charts of $X$ are monomial maps, there is a sense to asking for the points of $X$ with coordinates in 
$\RR_{\geq 0} := \left\{x\in \RR \st x\geq 0\right\}$, a fact first realized by Ehlers \cite[\S IV]{E}.

Intrinsically, for each cone $\sigma\in \Delta_{X}$ with associated semigroup $S_{\sigma}:=\sigma^{\vee}\cap M$,
the $\RR_{\geq 0}$-points of the affine variety $U_{\sigma}\subseteq X$ may be identified with the semigroup
homomorphisms $\Hom_{\sg}(S_{\sigma},\RR_{\geq 0})$, where $\RR_{\geq 0}$ is thought of as a semigroup under 
multiplication. 
One may similarly look at the points of $X$ with coordinates in  $\RR_{>0} := \left\{x\in \RR \st x>0\right\}$.

We use $X_{\geq 0}$ and $X_{>0}$ respectively for the corresponding subspaces of $X(\CC)$.

As elementary examples, we note that 
$$
\PP^{n}_{\geq 0} = \left\{ [x_0:x_1:\cdots : x_n] \st x_i\in \RR_{\geq 0} \right\}
= \left\{ [x_0:x_1:\cdots : x_n] \st x_i\in \RR_{\geq 0},\,\sum x_i = 1 \right\}
$$
is the $n$-dimensional simplex, and that when $X=(\PP^1)^{n}$,  $X_{\geq 0} = 
(\PP^1_{\geq 0})^n \cong [0,1]^n$ is the $n$-dimensional cube. 

\point 
\label{sec:intro-description}
The decomposition $\CC^{*} = S^1\times \RR_{>0}$, where $S^1$ denotes the unit circle, induces a 
similar decomposition of the torus $T_{X}$ of $X$ into its compact and non-compact parts~:
$$T_{X}\cong (S^1)^n \times (\RR_{>0})^n.$$

As shown in \cite[Chap.\ I]{AMRT} the exponential map gives a 
natural identification between $N_{\RR}$ and the non-compact part of $T_{X}$ above.

Via this exponential map, and the action of $T_{X}$ on $X$, $N_{\RR}$ acts on $X_{\geq 0}$, breaking it into
finitely many orbits, each one corresponding to a cone of $\Delta_{X}$ in a way parallel to the usual decomposition
of $X$ under the action of $T_{X}$ (see \cite[loc. cit.]{AMRT}). 

Specifically, for each cone $\sigma\in \Delta_{X}$, the stabilizer of a point in the orbit corresponding to $\sigma$
is the subspace $\langle \sigma\rangle\subseteq N_{\RR}$, and so the orbit is isomorphic to 
$N_{\RR}/\langle \sigma\rangle$.  Thus $X_{\geq 0}$ decomposes naturally into open cells with the largest cell,
$X_{>0}$, identified with $N_{\RR}$ itself. 

Now assume in addition that $X$ is proper.
Our goal in this note is twofold~:

\EulerList
\begin{enumerate}
\item To show that $X_{\geq 0}$ is homeomorphic to the $n$-dimensional ball.

\smallskip
\item To show that these cells give $X_{\geq 0}$ the structure of a regular cell complex.
\end{enumerate}
\AlphaList

The key step is to show (\Euler{1}), which then easily implies (\Euler{2}).   Both of these
results are established in Theorem \ref{thm:top-structure-of-X} below.

When $X$ is a smooth and proper, but not necessarily projective, toric variety, the
fact that $X_{\geq 0}$ is homeomorphic to a ball was first shown by Ehlers, \cite[p. 143, Lemma 1]{E}.
The case of a projective, but possibly singular, toric variety was established by Jurkiewicz \cite{J}, see
also \cite[p. 81, Proposition]{F}.  In this case the topological space $X_{\geq 0}$ may be identified
with the image of the moment map, \cite[p.\ 82--83]{F}.

The result for the common generalization, that of a possibly singular proper toric variety, is 
presumably ``well-known'' (or at least ``well-expected''), but the author has been unable to find such an argument 
in the literature, which is why one is provided here. 

\point
The motivation for this result is for use in \cite{RZ}.  There, Theorem \ref{thm:top-structure-of-X} is a key input 
in the construction of a reduced \v{C}ech complex for any semi-proper toric variety (that is a toric variety $X$
such that the natural map $X\longrightarrow \Spec(\Gamma(X,\Osh_{X}))$ is proper). 
This complex has the properties that ({\em i}) all open sets in the complex are torus stable open affines; 
({\em ii}) that the complex is of the same dimension as the cohomological dimension of $X$; 
and ({\em iii}) that this complex computes the cohomology for all coherent sheaves on $X$ 
(see \cite[Theorem 4.16]{RZ}).
In \cite{RZ} this complex is used in turn to give an algorithm computing the higher direct images of line bundles
under toric fibrations between smooth proper toric varieties, where the target is also projective. 

As part of the proof of Theorem \ref{thm:top-structure-of-X} we study the partition of $N_{\RR}$ given by the 
barycentric decomposition of $\Delta_{X}$. 
This note would be slightly shorter if we used the pieces of the barycentric subdivision as the cells instead of
the $N_{\RR}$-orbits,  but this has the disadvantage of introducing many more cells.  
For example, for $X=(\PP^1)^n$, this decomposition would result in $2^n\cdot n! $ cells of dimension $n$, 
whereas the theorem gives a cell decomposition with only a single $n$-dimensional cell, namely $X_{>0}\cong N_{\RR}$.
This reduction in the number of cells is very useful in the performance of the algorithm, and in any case
the decomposition of $X_{\geq 0}$ into $N_{\RR}$-orbits seems the most natural way to give $X_{\geq 0}$ the
structure of a cell complex.

Let $B^{n}$ be the compactification of $N_{\RR}$ by the sphere at infinity (see \S\ref{sec:definition-of-BN}).
Besides establishing (\Euler{1}) and (\Euler{2}) above, two other results in this note which may be of interest are~:

\begin{itemize}
\item The construction of a homeomorphism $N_{\RR}\cong X_{>0}$ which extends to a homeomorphism 
$B^{n}\cong X_{\geq 0}$. (Note that the homeomorphism
$N_{\RR}\cong X_{>0}$ given by the exponential map mentioned in \S\ref{sec:intro-description}
does not extend to a map from $B^{n}$ to $X_{\geq 0}$, see 
\S\ref{sec:Remarks-on-extension-of-homeomorphism}(\Euler{1}).)  
This homeomorphism is part of the statement of Theorem \ref{thm:top-structure-of-X}. 

\medskip
\item For each $n$-simplex in the barycentric decomposition of $B^{n}$ induced by $\Delta_{X}$, an explicit 
parameterization of the corresponding simplex in $X_{\geq 0}$. Combining these, one obtains a 
decomposition of $X_{\geq 0}$ into parameterized simplices.
See \S\ref{sec:remark-on-explicit-parameterization} and \S\ref{sec:Remarks-on-extension-of-homeomorphism}(\Euler{2}).
\end{itemize}

For a heuristic explanation of why one might try and analyze this problem by studying the barycentric subdivision
see \S\ref{sec:Remarks-on-extension-of-homeomorphism}(\Euler{3}).

\bpoint{Organization of the paper}
In \S\ref{sec:AMRT} we recall some essential results from \cite{AMRT} used repeatedly in the rest of the paper.
Sections \S\ref{sec:barycentric-subdivision}  and \S\ref{sec:definition-of-BN} recall standard
definitions and also serve to establish the notation used.  The rest of the paper is a sequence of constructions
and intermediate results leading up to the proof of Theorem \ref{thm:top-structure-of-X}. 
The most involved of these are Theorems \ref{thm:C-F-monomial-diagram} and \ref{thm:intersection-of-cones}, 
studying the cones of the barycentric subdivision of $N_{\RR}$, their closures in $X_{\geq 0}$, and their
intersections.

\bpoint{Acknowledgements}
The impetus for this result originated in the paper \cite{RZ}, and I thank my coauthor Sasha Zotine for 
pointing out that for a projective toric variety the polytopes associated respectively to $\Delta_{X}$ and to the 
torus orbits are dual, a remark which led to the strategy of the proof.

\section{Background material and proof of the main theorem} 

\bpoint{Results from ``Smooth compactifications of locally symmetric varieties''} 
\label{sec:AMRT}
Let $X$ be a complex toric variety, with torus $T_{X}$ and fan $\Delta_{X}\subset N_{\RR}$.  For an element 
$\alpha\in M$ we use $\charx^{\alpha}\colon T_{X}\longrightarrow \CC^{*}$ for the corresponding torus character.
We use $\langle \smallbullet,\smallbullet\rangle \colon M\times N \longrightarrow \ZZ$ for the pairing between
$M$ and $N$, or its extensions after tensoring with a field. 

The following results are taken directly from \cite[Chap. I]{AMRT}. 

\EulerList
\begin{enumerate}
\item There is a unique homomorphism $\exp\colon N_{\CC}\longrightarrow T_{X}$ 
such that for all $\alpha\in M$, and $z\in N_{\CC}$, $\charx^{\alpha}(\exp(z)) = e^{2\pi i \langle \alpha,z\rangle}$.
\PauseEnumerate

The kernel of $\exp$ is $N\subset N_{\CC}$ and identifies $N$ with the fundamental group of $T_{X}$.

\ResumeEnumerate
\item Setting $n=\dim(X)$ and writing $N_{\CC} = N_{\RR} \oplus iN_{\RR}$, $\exp$ thus gives an isomorphism
$$T_{X} \cong \left({N_{\CC}}/{N}\right) \oplus i N_{\RR} \cong (S^1)^n\oplus (\RR_{> 0})^n.$$
The two terms are the compact and non-compact parts of $T_{X}$, respectively. 
We use $T_{X,>0}$ for the non-compact part of $T_{X}$, a notation justified by (\Euler{3}) below. 

\PauseEnumerate

Let $\expi\colon N_{\RR} \longrightarrow T_X$ denote\footnote{ This name does not appear in \cite{AMRT}, and
is being introduced by the author. The author does not think that it is a very good name, but 
did not find a better one.} the composition of multiplication by $i$ on $N_{\RR}$ 
and the map $\exp$ above.  Then $\expi$ gives an isomorphism between $N_{\RR}$ and the non-compact part
of $T_{X}$.  

\ResumeEnumerate
\item  For each cone $\sigma\in \Delta_{X}$, $U_{\sigma,>0}$ is the image of 
$N_{\RR}\stackrel{\expi}{\cong} T_{X,>0}$ 
under the open embedding
$T_{X}\hookrightarrow U_{\sigma}$. 
\PauseEnumerate

Thus, since the $U_{\sigma}$ cover $X$, $N_{\RR}\stackrel{\expi}{\cong} X_{>0}$, a fact we will use frequently.

\ResumeEnumerate
\item $X_{\geq 0}$ is the topological closure of $X_{>0}$ in $X$.
\PauseEnumerate

Point (\Euler{4}) is not explicitly stated in \cite[Chap. I]{AMRT}, but follows immediately from the description
of the $N_{\RR}$-orbits (acting via $\expi$) on the points of $U_{\sigma,\geq 0}$ for each 
cone $\sigma$. 

We will also frequently need the description of the image of a point of $N_{\RR}$ under $\expi$ followed
by a closed embedding.

Let $\sigma\in \Delta_{X}$ be a cone, and $S_{\sigma}=\sigma^{\vee}\cap M$ the corresponding semigroup.
Choosing generators $\alpha_1$, \ldots, $\alpha_m$ for $S_{\sigma}$, the surjection
$\CC[y_1,\ldots, y_m]\longrightarrow \CC[S_{\sigma}]$ given by $y_i\mapsto \charx^{\alpha_i}$ gives a closed
immersion of $U_{\sigma}:=\Spec(\CC[S_{\sigma}])$ into $\AA^{m}_{\CC}$.  

For any $x\in N_{\RR}$, the image of $\expi(x)$ under this closed immersion is the point
$$\left(e^{-2\pi \langle \alpha_1,x\rangle}, e^{-2\pi \langle \alpha_2,x\rangle},\,\ldots,\,
e^{-2\pi \langle \alpha_m,x\rangle}\right).$$

Finally, we note that, since $X_{\geq 0}$ is closed in $X$, if $X$ is proper, then $X_{\geq 0}$ is compact.

\AlphaList

\bpoint{Barycentric subdivisions}
\label{sec:barycentric-subdivision} 
Given the fan $\Delta_{X}$, the associated barycentric subdivision is a subdivision of $N_{\RR}$ 
into simplicial cones which is a refinement of 
the decomposition into cones given by $\Delta_{X}$.  For each cone $\sigma\in \Delta_{X}$ we put
$$\Bsub_{\sigma} := \sum_{\rho\in \sigma(1)} v_{\rho},$$
where $v_{\rho}$ denotes the integral generator of the ray $\rho$, and $\sigma(1)$ the rays of $\sigma$. 

We note that the notation is unambiguous when $\sigma$ is itself a ray.  Writing instead (temporarily)  
$\Bsub'_{\sigma} = \sum_{\rho\in \sigma(1)} v_{\rho}$,  if $\sigma$ is a ray, $\sigma=\rho$, 
then $\Bsub'_{\rho}=v_{\rho}$.

For each flag $\Flag = \left\{ \sigma_{i_1} \subset \sigma_{i_2} \subset \cdots \subset \sigma_{i_k}\right\}$ of
cones in $\Delta_{X}$ we denote by $C_{\Flag}$ the $k$-dimensional simplicial cone spanned by 
$\Bsub_{\sigma_{i_1}}$, \ldots, $\Bsub_{\sigma_{i_k}}$, i.e., the subset 
$$C_{\Flag} := \left\{ \sum_{j=1}^{k} u_j \Bsub_{\sigma_{i_j}} \st u_j\in \RR_{\geq 0}\right\} \subset N_{\RR}.$$
In our notation for flags we do not include the cone $\{0\}$, and the inclusion 
$\sigma_{i_{j}}\subset \sigma_{i_{j+1}}$ must be strict.  We do use the empty flag, $\Flag=\emptyset$, and in this case
set $C_{\emptyset} = \{0\} \subset N_{\RR}$. 

Given two flags $\Flag_1$ and $\Flag_2$, we denote by $\Flag_1\cap \Flag_2$ the flag consisting of the cones
common to both $\Flag_1$ and $\Flag_2$, including the possibility that $\Flag_1\cap \Flag_2=\emptyset$.   
From these definitions, and the definition of a fan, one can show that 
$C_{\Flag_1}\cap C_{\Flag_2} = C_{\Flag_1\cap \Flag_2}$, and
that for any cone $\sigma\in \Delta_{X}$, 
$$|\sigma| = \bigcup_{\Flag} C_{\Flag},$$
where the union is over flags $\Flag$ contained in $\sigma$, and $|\sigma|$ denotes the support of $\sigma$. 
Setting $k=\dim(\sigma)$, a minimal union is obtained by 
taking the flags $\Flag=\{\sigma_{i_1}\subset \sigma_{i_2} \subset \cdots \subset \sigma_{i_{k-1}} \subset \sigma\}$, 
i.e., the  union of the $k$-dimensional cones whose top element in each flag is $\sigma$. 

We will use the decomposition of $N_{\RR}$ into the cones $C_{\Flag}$ to understand the 
compactification $N_{\RR}\hookrightarrow X_{\geq 0}$.  We note that we are not changing the fan $\Delta_{X}$, i.e.,
we are not blowing up $X$, but rather studying the compactification of $N_{\RR}$ given by $\Delta_{X}$ (i.e., in
$X_{\geq 0}$), ``cone by cone''.

\bpoint{Definition of $\theta$} 
\label{sec:def-of-theta}
Fix $n\geq 1$, and for each $j$, $1\leq j\leq n$, define $\theta_j\colon [0,1]^n\longrightarrow [0,1]$ by
$\theta_{j}(z_1,\ldots, z_n) = \prod_{i=j}^{n} z_i$. 
We then define $\theta\colon [0,1]^n\longrightarrow [0,1]^n$ to be the product map $\theta = \prod_{j=1}^{n} \theta_j$.
For instance, when $n=4$ this map is 
$$(z_1,z_2,z_3,z_4) \xmapsto{\rule{0.25cm}{0cm}\theta\rule{0.25cm}{0cm}} (z_1z_2z_3z_4,\, z_2z_3z_4,\,z_3z_4,\,z_4).$$

We use $w_1$,\ldots, $w_n$ for the coordinates on the target copy of $[0,1]^n$, and (as above) $z_1$,\ldots, $z_n$
for the coordinates on the source copy of $[0,1]^n$. 

\bpoint{Proposition} 
\label{prop:properties-of-theta}
Fix $n\geq 1$ and let $\theta\colon [0,1]^{n}\longrightarrow [0,1]^{n}$ be the map of \S\ref{sec:def-of-theta}. 
Then

\begin{enumerate}
\item The image of $\theta$ is the $n$-dimensional simplex 
$$\Delta_n:=\left\{(w_1,\ldots, w_n) \st 0\leq w_1\leq w_2\leq \cdots \leq w_{n-1}\leq w_{n} \leq 1
\rule{0cm}{0.35cm}\right\}\subset [0,1]^n.$$

\smallskip
\item The composition $(0,1]^n\hookrightarrow [0,1]^n\stackrel{\theta}{\longrightarrow} [0,1]^n$ induces a homeomorphism
of $(0,1]^n$ with $\Delta_n\cap \{(w_1,\ldots, w_n)\st w_1\neq 0\}$. 

\smallskip
\item Let $\{b_{ij}\}$, $i$, $j\in \{1,\ldots, n\}$ be non-negative integers, with $b_{ii}>0$ and $b_{ij}=0$ if 
$j<i$.  Define a monomial map $\psi\colon [0,1]^n\longrightarrow [0,1]^n$ by setting the $i$-th coordinate
of the map to be $(w_1,\ldots, w_n)\mapsto \prod_{j=1}^{n} w_j^{b_{ij}}$.  
I.e, the map is
$$(w_1,\ldots, w_n) \mapsto \left(w_1^{b_{11}}w_2^{b_{12}}\cdots w_{n}^{b_{1n}},\,\ldots,\, 
w_{n-1}^{b_{n-1,n-1}}w_{n}^{b_{n-1,n}},\, w_n^{b_{nn}}\right).$$
Then $\psi|_{\Delta_n}$ is injective.
\end{enumerate}

\proof 
The arguments for ({\em a}) and ({\em b}) are straightforward, and we leave them to the reader.
To see ({\em c}), suppose that $(w_1,\ldots, w_n)\in \Delta_n$.

If none of the entries of $\psi(w_1,\ldots, w_n)$ are equal to zero, then none of the $w_i$ are equal to zero,
and we can recover the $w_i$ by using the upper triangular structure of the map and extracting roots. (E.g.,
we recover $w_n$ by taking $b_{nn}$-th root of the last entry, and then recover $w_{n-1}$ by dividing
the second last entry by $w_{n}^{b_{n-1,n}}$ and then taking the $b_{n-1,n-1}$-st root.) 

Otherwise, let $i$ be the largest index such that the $i$-th entry of $\psi(w_1,\ldots, w_n)$ is zero. 
From the upper triangular structure of the map, and the fact that no higher entry is zero, we conclude
that $w_i=0$ and $w_j\neq 0$ for $j>i$.  
We may then apply the previous recursive procedure to determine $w_j$ from the nonzero 
entries of $\psi(w_1,\ldots, w_n)$, for those $j>i$. 
On the other hand, by the inequalities in ({\em a}) defining $\Delta_n$,  if $w_i=0$ then $w_j=0$
for all $j\leq i$.   Thus, for $(w_1,\ldots, w_n)\in \Delta_n$,  
knowing $\psi(w_1,\ldots, w_n)$ is sufficient
to determine $(w_1,\ldots, w_n)$.  I.e., $\psi|_{\Delta_n}$ is injective. \epf

\point
For later use, we record the formula for points of $\Delta_n$ in barycentric coordinates.   
For $i$ such that
$0\leq i\leq n$, set $p_i:=(w_1,\ldots, w_n)$, with $w_j=0$ when $j\leq i$ and $w_j=1$ when $j>i$.  Thus
$p_0=(1,1,\ldots, 1)$, $p_1=(0,1,1,\ldots, 1)$, up to $p_{n-1}=(0,0,\ldots, 0,1)$ and $p_n=(0,0,\ldots, 0)$.

The simplex $\Delta_n$ is the convex hull of $p_0$,\ldots, $p_n$. 
Given $\xi_0$, $\xi_1$, \ldots, $\xi_n$, with each $\xi_i\geq 0$ and $\sum_{i=0}^{n} \xi_i = 1$, the corresponding
convex combination of $p_0$,\ldots, $p_n$ is 

\begin{equation}
\label{eqn:simplex-coordinates}
\sum_{i=0}^{n} \xi_i p_i = \left(\xi_0,\, \xi_0+\xi_1,\, \xi_0+\xi_1+\xi_2,\, \ldots,\,
\xi_0+\xi_1+\cdots + \xi_{n-1} \rule{0cm}{0.35cm}\right),
\end{equation}

i.e., the point whose $j$-th coordinate is $\sum_{i=0}^{j-1}\xi_i$.
Conversely, the barycentric coordinates for a point $(w_1,\ldots, w_n)\in \Delta_n$, are given by
$\xi_0=w_1$, $\xi_j=w_{j+1}-w_{j}$ for $j=1$,\ldots, $n-1$, and $\xi_{n} = 1-w_n$.

\bpoint{Theorem} 
\label{thm:C-F-monomial-diagram}
Let $\Flag = \{\sigma_1 \subset \sigma_2\subset \cdots \subset \sigma_n\}$ be a maximal flag in $\Delta_{X}$,
and $C_{\Flag}\subset N_{\RR}$ the corresponding simplicial cone.  
As in \S\ref{sec:AMRT}, let $\alpha_1$, \ldots, $\alpha_m$ be generators of $\sigma_n^{\vee}\cap M$,
and $U_{\sigma_n}\hookrightarrow \AA^{m}$ the corresponding closed embedding.
Let $\theta$ be the map defined in \S\ref{sec:def-of-theta} and, as above, set $\Delta_n:=\Im(\theta)$. 

Let $\beta_1$,\ldots, $\beta_n\in N_{\QQ}$ be the basis dual to $\Bsub_{\sigma_1}$, \ldots, $\Bsub_{\sigma_n} \in M$.
Define the homeomorphism $\exp_{\Flag}\colon C_{\Flag} \stackrel{\sim}{\longrightarrow} (0,1]^n$  by 

\begin{equation}
\label{eqn:exp-F-formula-1}
\exp_{\Flag}(x) = \left(e^{-2\pi\langle \beta_1,x\rangle},\, e^{-2\pi\langle \beta_2,x\rangle},\,\ldots,\,
e^{-2\pi\langle \beta_n,x\rangle}\right).
\end{equation}

I.e., writing each $x\in C_{\Flag}$ in simplicial coordinates on $C_{\Flag}$, $x=\sum_{i=1}^{n} u_i\Bsub_{\sigma_i}$,
the map 

\begin{equation}
\label{eqn:exp-F-formula-2}
(u_1,\ldots, u_n) \mapsto \left(e^{-2\pi u_1},\, e^{-2\pi u_2},\,\ldots, e^{-2\pi u_n}\right).
\end{equation} 

Then there exists a monomial map $\psi\colon [0,1]^n\longrightarrow \AA^m$ such that the diagram  

\begin{equation}
\label{eqn:C-F-closure-factorization}
\rule{3cm}{0cm}
\xymatrix{
C_{\Flag}\rule[-1.5mm]{0cm}{1.5mm}\ar[rr]_{\sim}^{\exp_{\Flag}}\ar@{^{(}->}[d]_{\mbox{{\tiny closed}}} & & (0,1]^n\ar@{^{(}->}[dd]^{h} \\
N_{\RR}\rule[-1.5mm]{0cm}{1.5mm}\ar@{^{(}->}[d]^{\expi}_{\mbox{\tiny open}}  & & \\
U_{\sigma_n}\rule[-1.5mm]{0cm}{1.5mm}\ar@{^{(}->}[d]_{\mbox{\tiny closed}} & & [0,1]^n\ar[d]^{\theta} &  & \mbox{\small (coordinates $z_1$, \ldots,$z_n$)}\\
\AA^{m} & & [0,1]^n\ar[ll]_{\psi} & & \mbox{\small (coordinates $w_1$, \ldots, $w_n$)} \\
}
\end{equation}

commutes. Moreover, $\Im(\psi\circ\theta)$ is contained in $U_{\sigma_n}$, and $\psi|_{\Delta_{n}}$ induces a 
homeomorphism of $\Delta_{n}$  with $\overline{C}_{\Flag}$, where $\overline{C}_{\Flag}$ 
denotes the topological closure of $C_{\Flag}$ in $U_{\sigma_n}$.

\bpf
Let $f\colon C_{\Flag}\longrightarrow \AA^{m}$ be the composite of the left-hand vertical morphisms.
From the definition of $\expi$, and the closed
immersion $U_{\sigma_{n}}\hookrightarrow \AA^{m}$ it follows that for each $x\in C_{\Flag}$, 
$$f(x) = \left(e^{-2\pi \langle \alpha_1,x\rangle},\ldots, e^{-2\pi \langle \alpha_m,x\rangle}\right).$$
Equivalently, using $y_1$,\ldots, $y_m$ for the coordinates on $\AA^{m}$, for each $i\in \{1,\, \ldots,\, m\}$, 

\begin{equation}
\label{eqn:f-yi-pullback}
(f^{*}(y_i))(x) = e^{-2\pi \langle \alpha_i,x\rangle}.
\end{equation}

Let $h$ denote the open immersion $(0,1]^n\hookrightarrow [0,1]^n$.  Then for each $j\in \{1,\ldots, n\}$,
$$((h\circ \exp_{\Flag})^{*}z_j)(x) = e^{-2\pi \langle \beta_j,x\rangle}.$$
From the definition of $\theta$ we then get that 

\begin{equation}
\label{eqn:wj-pullback}
((\theta\circ h\circ \exp_{\Flag})^{*}w_j)(x) = e^{-2\pi\langle \sum_{k=j}^{n} \beta_k,x\rangle}.
\end{equation} 

Saying that $\psi$ is a monomial map is equivalent to saying that, for each $i$, 
$\psi^{*} y_i = w_1^{b_{i1}}w_2^{b_{i2}}\cdots w_{n}^{b_{in}}$ for some nonnegative integers $b_{i1}$, \ldots, $b_{in}$.

Comparing \eqref{eqn:f-yi-pullback} and \eqref{eqn:wj-pullback}, the condition on the $b_{ij}$ to make the 
diagram commute is that for each $i$, 

\begin{equation}
\label{eqn:alpha-commutativity-condition}
\alpha_i = \sum_{j=1}^{n} b_{ij}\sum_{k=j}^{n} \beta_k 
= \sum_{k=1}^{n} \left(\sum_{j=1}^{k} b_{ij}\right) \beta_k.
\end{equation}

Since $\beta_1$,\ldots, $\beta_n$ are a basis for $M_{\QQ}$, for each $i\in \{1,\ldots, m\}$  we may write 
$\alpha_{i} = \sum_{k=1}^{n} c_{ik} \beta_k$ for unique $c_{ik}\in \QQ$. 
Moreover, since $\beta_1$,\ldots, $\beta_n$ are dual to $\Bsub_{\sigma_1}$,\ldots, $\Bsub_{\sigma_n}$, 
$c_{ik} = \langle \alpha_i,\Bsub_{\sigma_k}\rangle$. 
From this formula for the $c_{ik}$ we draw two conclusions~:

\RomanList
\begin{enumerate}
\item Each $c_{ik}$ is a nonnegative integer.
\item For each $k\in \{2,\ldots, n\}$, $c_{ik}\geq c_{i,k-1}$.
\end{enumerate}
\AlphaList

Conclusion ({\em i}) follows since $\Bsub_{\sigma_k}\in N\cap \sigma_n$ and $\alpha_i\in \sigma_n^{\vee}\cap M$. 
To see ({\em ii}) we recall that each $\Bsub_{\sigma_k}$ is the sum 
$\Bsub_{\sigma_k} = \sum_{\rho\in \sigma_k(1)} v_{\rho}$, i.e., the sum of the generators of the rays 
of $\sigma_k$.  Since $\sigma_{k-1}$ is a face of $\sigma_k$, it follows
that $\Bsub_{\sigma_{k}}-\Bsub_{\sigma_{k-1}}$ is a sum of the generators of those rays of $\sigma_{k}$ which are not
rays of $\sigma_{k-1}$.  In particular, $\Bsub_{\sigma_{k}}-\Bsub_{\sigma_{k-1}}$ is a sum, 
with non-negative coefficients, of elements of $\sigma_n$.  Since $\alpha_i\in\sigma_n^{\vee}$, we therefore have
$\langle \alpha_i,\Bsub_{\sigma_k}-\Bsub_{\sigma_{k-1}}\rangle\geq 0$, i.e., $c_{ik}\geq c_{i,k-1}$. 

Setting $b_{i1}=c_{i1}$, and, for $j\geq 2$, $b_{ij} = c_{ij}-c_{i,j-1}$, we have by ({\em i}) and ({\em ii}) above
that each $b_{ij}$ is a non-negative integer. From the definition of the $b_{ij}$ we have 
$c_{ik}=\sum_{j=1}^{k} b_{ij}$, and thus 
\eqref{eqn:alpha-commutativity-condition} holds.  Thus the $b_{ij}$ define a monomial map $\psi$ which makes
\eqref{eqn:C-F-closure-factorization} commute.   The proof moreover shows that the $b_{ij}$ are unique, and thus that
$\psi$ is unique. 

We now turn to the proof of the second part of the theorem.

The fact that $\Im(\psi\circ\theta)$ is contained in $U_{\sigma_n}$ follows from the commutativity of the
diagram and the fact that $U_{\sigma_n}$ is closed in $\AA^{m}$.  Specifically, by the commutativity of the
diagram, $(\psi\circ\theta)|_{\Im(h)}$ induces a homeomorphism of $(0,1]^n$ with $C_{\Flag}\subseteq \AA^{m}$.
Since $(0,1]^n$ is dense in $[0,1]^n$, and since $[0,1]^n$ is compact, 
we conclude that $\Im(\psi\circ\theta)$ is the closure of $C_{\Flag}$
in $\AA^{m}$.  But, since $U_{\sigma_n}$ is closed in $\AA^{m}$, this closure is $\overline{C}_{\Flag}$. 
This shows both that $\Im(\psi\circ\theta)$ is contained in $U_{\sigma_n}$, and that $\Im(\psi\circ\theta)=
\overline{C}_{\Flag}$.

Thus, $\psi|_{\Delta_n}$ is a surjective map from $\Delta_{n}$ to $\overline{C}_{\Flag}$.  Since 
$\Delta_{n}$ is compact and $\overline{C}_{\Flag}$ Hausdorff ($\overline{C}_{\Flag}$ is closed in $\AA^{m}$), 
to show that $\psi|_{\Delta_{n}}$ is a homeomorphism onto its image it suffices to show that $\psi|_{\Delta_{n}}$
is injective, and to do this we reduce to a situation where we may apply
Proposition \ref{prop:properties-of-theta}({\em c}).  

Given $n$ coordinates $y_{i_1}$,\ldots, $y_{i_n}$ we may
project onto these coordinates to define a map $\pi\colon \AA^{m}\longrightarrow \AA^{n}$, and obtain a monomial map
$\psi':=\pi\circ\psi$.  If $\psi'$ satisfies the conditions of 
Proposition \ref{prop:properties-of-theta}({\em c}), then by the proposition $\psi'$, and hence $\psi$, is injective
on $\Delta_{n}$.   Thus we are reduced to finding such a set of coordinates.

Since $\Im(\psi\circ\theta)$ lies in $U_{\sigma_n}$, the claim on injectivity is independent of the embedding
$U_{\sigma_n}\hookrightarrow \AA^{m}$, and so independent of the generators chosen for $\sigma^{\vee}\cap M$.
In particular, we may add extra elements to the generating set without changing the truth or 
falsity of the injectivity statement.  We may also clearly reorder the generators.

The faces $\sigma_1\subset \sigma_2 \subset \cdots \subset \sigma_{n-1}$ of $\sigma_{n}$ give corresponding
faces of the dual cone 
$$\sigma'_{n-1} \subset \sigma'_{n-2} \subset \cdots \subset \sigma'_{1} \subset \sigma_n^{\vee},$$
where $\sigma'_{i} := \left\{\alpha\in \sigma_{n}^{\vee} \st \alpha|_{\sigma_i}\equiv 0\right\}$.

Choose $\alpha_n$ in the relative interior of $\sigma'_{n-1}$, $\alpha_{n-1}$ in the relative interior of
$\sigma'_{n-2}$, down to $\alpha_2$ in the relative interior of $\sigma'_{1}$, and $\alpha_1$ in the interior
of $\sigma^{\vee}$, with all choices also in $M$.   With these choices, $\alpha_i|_{\sigma_j}\equiv 0$ for 
all $j<i$, and $\alpha_i|_{\sigma_i}\not\equiv 0$.  
From this and the formula $c_{ij}=\langle \alpha_i,\Bsub_{\sigma_j}\rangle$ we conclude that $c_{ij} =0$ for all $j<i$, 
and that $c_{ii}> 0$.  
Thus, with $b_{i1}=c_{i1}$, and $b_{ij}=c_{ij}-c_{i,j-1}$ for $j\geq 2$, the $b_{ij}$
satisfy the conditions of Proposition \ref{prop:properties-of-theta}({\em c}).  
Including the $\alpha_i$, $i=1$, \ldots, $n$, among the generators of $\sigma^{\vee}\cap M$, and 
setting $\pi$ to be the projection onto the corresponding coordinates, we conclude that $\psi$ is injective on the 
image of $\theta$. 
\epf

\bpoint{Corollary} 
\label{cor:structure-of-C-F}
Let $\Flag = \{\sigma_1\subset \cdots \subset \sigma_n\}$ be a maximal flag in $\Delta_{X}$, 
and $\overline{C}_{\Flag}$ the topological closure of $C_{\Flag}$ in $U_{\sigma_n}$, 
as in Theorem \ref{thm:C-F-monomial-diagram} above. 

\begin{enumerate}
\item $\overline{C}_{\Flag}$ is an $n$-dimensional simplex.
\item $\overline{C}_{\Flag}$ is closed in $X$.
\item The part of $\overline{C}_{\Flag}$ in the boundary of $X_{\geq 0}$, i.e., 
$$\overline{C}_{\Flag}\cap (X_{\geq 0}\setminus X_{>0}) = \overline{C}_{\Flag}\setminus C_{\Flag},$$
is an $(n-1)$-dimensional simplex.
\item $X_{\geq 0} = \displaystyle\bigcup_{\mbox{\scriptsize $\Flag$ maximal}} \overline{C}_{\Flag}$.
\end{enumerate}

\bpf
({\em a}) By Theorem \ref{thm:C-F-monomial-diagram} the map $\psi$ (constructed relative to this flag $\Flag$)
gives a homeomorphism of $\Delta_{n}$ with $\overline{C}_{\Flag}$.  By 
Proposition \ref{prop:properties-of-theta}({\em a}) $\Delta_{n}$ is an $n$-dimensional simplex. 

Since $\overline{C}_{\Flag}$ is compact, its image under the composition of the inclusions 
$\overline{C}_{\Flag}\hookrightarrow U_{\sigma_n} \hookrightarrow X$ is closed in $X$,  giving ({\em b}). 

We first justify the claimed equality in ({\em c}). Since $C_{\Flag}$ is closed in $N_{\RR}$, and 
$N_{\RR}\cong X_{>0}$, $\overline{C}_{\Flag}\setminus C_{\Flag}$ lies in $X_{\geq 0}\setminus X_{>0}$, and thus 
$$\overline{C}_{\Flag}\cap (X_{\geq 0}\setminus X_{>0}) = \overline{C}_{\Flag}\setminus C_{\Flag}.$$

Again using the homeomorphism $\psi|_{\Delta_n}$ between $\Delta_{n}$ and $\overline{C}_{\Flag}$,
by Proposition \ref{prop:properties-of-theta}({\em a}) and ({\em b}), $\overline{C}_{\Flag}\setminus C_{\Flag}$
may be identified with the set 
$$\left\{(w_1,\ldots, w_n) \in [0,1]^n \st 0 =  w_1\leq w_2\leq \cdots \leq w_{n-1}\leq w_{n} \leq 1
\rule{0cm}{0.35cm}\right\},$$
giving ({\em c}). 

Recall (\S\ref{sec:AMRT}) that $X_{\geq 0}$ is the closure of $X_{>0}\cong N_{\RR}$ in $X$. 
Since $N_{\RR} = \bigcup_{\mbox{\scriptsize $\Flag$ maximal}} C_{\Flag}$, ({\em d}) follows immediately. \epf

\bpoint{Remark} 
\label{sec:remark-on-explicit-parameterization} 
We note that as a by-product of the proof of Theorem \ref{thm:C-F-monomial-diagram} we obtain an explicit 
parameterization of each $n$-simplex $\overline{C}_{\Flag}\subset X_{\geq 0}$ for maximal flags $\Flag$.

Specifically, suppose we are given a maximal flag 
$\Flag = \{\sigma_{1} \subset \sigma_{2} \subset \cdots \subset \sigma_{n}\}$, and embedding 
$U_{\sigma_n}\hookrightarrow \AA^{m}$ given by a choice of generators $\alpha_1$,\ldots, $\alpha_m$ of
$\sigma_{n}^{\vee}\cap M$.  Setting 
$$
b_{ij} = 
\left\{
\begin{array}{cl}
\langle \alpha_i, v_{\sigma_1}\rangle & \mbox{if $j=1$, and} \\
\langle \alpha_i, v_{\sigma_{j}}-v_{\sigma_{j-1}}\rangle & \mbox{for $j=2,\ldots, n$} \\
\end{array}\right.,
$$
we then define a monomial map $\psi\colon [0,1]^n\longrightarrow \AA^{m}$ by
$$\psi(w_1,\ldots, w_n) = \left(\prod_{j=1}^{n} w_j^{b_{1j}},\, \prod_{j=1}^{n} w_{j}^{b_{2j}},\, \ldots,\,
\prod_{j=1}^{n} w_{j}^{b_{mj}}\right).$$
The proof of Theorem \ref{thm:C-F-monomial-diagram} shows that $\psi$ gives a homeomorphism of the $n$-simplex
$\Delta_n\subset [0,1]^n$ 
(as in Proposition \ref{prop:properties-of-theta}({\em a})) 
with $\overline{C}_{\Flag}\subseteq U_{\sigma_n}\subseteq X_{\geq 0}$.

\bpoint{Definition of $\overline{C}_{\Flag}$ for any flag} 
Let $\Flag$ be a flag in $\Delta_{X}$ (including possibly $\Flag=\emptyset$). We define $\overline{C}_{\Flag}$
to be the topological closure of $C_{\Flag}$ in $X_{\geq 0}$, where we consider $C_{\Flag}$ a subset of $X_{>0}$
via the inclusion $C_{\Flag}\hookrightarrow N_{\RR}$ and the homeomorphism $\expi\colon N_{\RR} \cong X_{>0}$. 

For $\Flag$ a maximal flag, Corollary \ref{cor:structure-of-C-F}({\em b}) shows that this definition agrees with the 
definition of $\overline{C}_{\Flag}$ given in Theorem \ref{thm:C-F-monomial-diagram}.

If $\Flag$ is a maximal flag, and $\Flag_1\subset \Flag$ a subflag, then $C_{\Flag_1}$ is a face of $C_{\Flag}$,
and it follows from the explicit parameterization of $\overline{C}_{\Flag}$ that $\overline{C}_{\Flag_1}$ 
is a face of the $n$-simplex $\overline{C}_{\Flag}$ 
(e.g., see \eqref{eqn:simplex-boundary} in the proof of Theorem \ref{thm:intersection-of-cones}({\em a}) below). 
If $C_{\Flag_1}$ is a $k$-dimensional simplicial cone, then one also sees that $\overline{C}_{\Flag_1}$ is a 
$k$-simplex, and the part of $\overline{C}_{\Flag_1}$ in the boundary of $X_{\geq 0}$ is a $(k-1)$-simplex,
generalizing Corollary \ref{cor:structure-of-C-F}({\em c}) above.

Suppose that $\Flag_1$ and $\Flag_2$ are two flags in $\Delta_{X}$, and $\Flag_{12}:=\Flag_1\cap \Flag_2$ their
intersection.  Since 
$C_{\Flag_{12}}= C_{\Flag_1}\cap C_{\Flag_2}$, we have
$\overline{C}_{\Flag_{12}} \subseteq \overline{C}_{\Flag_1}\cap \overline{C}_{\Flag_2}$, however 
nothing so far shows that this is an equality.

For instance, to fix ideas, 
if $\Flag_{12}=\emptyset$, so that $C_{\Flag_1}$ and $C_{\Flag_2}$ intersect only
in the origin of $N_{\RR}$, so far we do not know that $\overline{C}_{\Flag_1}$ and $\overline{C}_{\Flag_2}$ cannot
have common points in 
$X_{\geq 0}\setminus X_{> 0}$. 

The arguments up until now are enough to show that the inclusion $N_{\RR}\hookrightarrow X_{\geq 0}$ extends to a map
from the $n$-dimensional ball onto $X_{\geq 0}$.  But, to know that this extension is injective on the boundary 
we need to establish equality in the intersection of cones above.

\bpoint{Theorem}
\label{thm:intersection-of-cones}

\begin{enumerate}
\item Let $\Flag_1$ be a maximal flag and $\Flag_2$ an arbitrary flag in $\Delta_{X}$, and put 
$\Flag_{12}:=\Flag_1\cap \Flag_2$.  Let $\{x_m\}_{m\geq 1}$ be a sequence of points of $C_{\Flag_2}$ which converges
to a point $x'$ of $\overline{C}_{\Flag_1}$.   Then $x'\in \overline{C}_{\Flag_{12}}$. 
\item The same as in ({\em a}) but with $\Flag_1$ an arbitrary flag in $\Delta_{X}$.
\item For any two flags $F_1$, $F_2$ in $\Delta_{X}$, 
$\overline{C}_{\Flag_1}\cap \overline{C}_{\Flag_2}=\overline{C}_{\Flag_1\cap \Flag_2}$. 
\end{enumerate}

\bpf The argument for ({\em a}) is the longest, but will quickly imply ({\em b}) and ({\em c}).
Suppose that $\Flag_1=\{\sigma_1\subset \cdots \subset \sigma_n\}$. 
As in the proof of Theorem \ref{thm:C-F-monomial-diagram}, let $\beta_1$,\ldots, $\beta_n\in M_{\QQ}$ be the 
basis dual to $\Bsub_{\sigma_1}$, \ldots, $\Bsub_{\sigma_n}$. 
Since $\Bsub_{\sigma_1}$, \ldots, $\Bsub_{\sigma_n}$ is a basis for $N_{\RR}$, 
each $x\in N_{\RR}$ has a unique expression $x=\sum_{i=1}^{n} u_i \Bsub_{\sigma_i}$, with each $u_i\in \RR$.
The formulas \eqref{eqn:exp-F-formula-1}  and \eqref{eqn:exp-F-formula-2} thus extend to give a homeomorphism
$N_{\RR}\stackrel{\sim}{\longrightarrow} (0,\infty)^n$ which we continue to call $\exp_{\Flag}$. 

We similarly promote the map $\theta\colon [0,1]^n\longrightarrow [0,1]^n$ of \S\ref{sec:def-of-theta} to
a map $[0,\infty)^n\longrightarrow [0,\infty)^n$, again called $\theta$, given by the same monomial formula.  

Let $\alpha_1$,\ldots, $\alpha_n\in \sigma_{n}^{\vee}$ be the elements constructed in the proof of 
Theorem \ref{thm:C-F-monomial-diagram} (i.e., those for which the associated $\{b_{ij}\}$ satisfy the conditions
of Proposition \ref{prop:properties-of-theta}({\em c})).  Extend these to generators
$\alpha_1$,\ldots, $\alpha_m$ of $\sigma_{n}^{\vee}\cap M$ and let $U_{\sigma_n}\hookrightarrow \AA^{m}$
be the associated closed embedding.  Let $\psi$ be the monomial map constructed in
Theorem \ref{thm:C-F-monomial-diagram} with respect to this embedding and the maximal flag $\Flag_{1}$. 
Using the same monomial formulas, we extend $\psi$ to a map $[0,\infty)^n\longrightarrow \AA^m$.

Combining these, we obtain the following diagram, extending \eqref{eqn:C-F-closure-factorization}.

\begin{equation}
\label{eqn:extended-factorization}
\rule{3.0cm}{0cm}
\xymatrix{
C_{\Flag}\rule[-1.5mm]{0cm}{1.5mm}\ar[rr]_{\sim}^{\exp_{\Flag}}\ar@{^{(}->}[d]_{\mbox{{\tiny closed}}} & & (0,1]^n\ar@{^{(}->}[d] \\
N_{\RR}\rule[-1.5mm]{0cm}{1.5mm}\ar@{^{(}->}[d]^{\expi}_{\mbox{\tiny open}}\ar[rr]_{\sim}^{\exp_{\Flag}}  & & (0,\infty)^n\ar@{^{(}->}[d]^{h} \\
U_{\sigma_n}\rule[-1.5mm]{0cm}{1.5mm}\ar@{^{(}->}[d]_{\mbox{\tiny closed}} & & [0,\infty)^n\ar[d]^{\theta}  \\
\AA^{m} & & [0,\infty)^n\ar[ll]_{\psi} & & \mbox{\small (coordinates $w_1$, \ldots, $w_n$)} \\
}
\end{equation}

The diagram \eqref{eqn:extended-factorization} again commutes, with the issue being to check that it commutes starting
at $N_{\RR}$.  Letting $y_1$,\ldots, $y_m$ be the coordinates on $\AA^{m}$, the pullback of $y_i$ to $N_{\RR}$ via
the vertical maps is, as before, the function $e^{-2\pi\langle \alpha_i,x\rangle}$, while the pullback of $y_i$
via $\psi\circ\theta\circ h\circ \exp_{\Flag}$ is, also as before, the function 
$e^{-2\pi \left\langle \sum_{k=1}^{n} \left(\sum_{j=1}^{k} b_{ij}\right) \beta_{k},x\right\rangle}$, 
where the $b_{ij}$ are the exponents defining $\psi$.   
However, these exponents have been chosen so that \eqref{eqn:alpha-commutativity-condition} holds, i.e., so that
the two pullbacks are equal.  Thus \eqref{eqn:extended-factorization} commutes.

Recall from Proposition \ref{prop:properties-of-theta}({\em a}) that $\Delta_n=\Im(\theta|_{[0,1]^n})$ is 
an $n$-simplex, given, as a subset of $[0,\infty)^n$, by the inequalities 

\begin{equation}
\label{eqn:simplex-boundary} 
0\stackrel{\infty}{\leq} w_1 \stackrel{u_1}{\leq} w_2 \stackrel{u_2}{\leq} \cdots 
\stackrel{u_{n-2}}{\leq} w_{n-1}\stackrel{u_{n-1}}{\leq} w_{n} \stackrel{u_n}{\leq} 1.
\end{equation}

The meaning of the symbols above the $\leq$ signs is the following. 
The faces of $C_{\Flag_1}$ are given by equations $u_j=0$, $j\in J$, for subsets $J\subseteq \{1,\ldots, n\}$. 
The faces of $\Delta_n$ are given by equalities in some subset of the conditions in \eqref{eqn:simplex-boundary}.  
The labels above match each $u_j$ to the appropriate inequality.  E.g., the face $u_1=0$ of $\overline{C}_{\Flag_1}$ 
becomes the face $w_1=w_2$ of $\Delta_n$.  The $(n-1)$-simplex ``at infinity'' of $\overline{C}_{\Flag_1}$ is
(via $\psi$) the face $0=w_1$ of $\Delta_n$ (Proposition \ref{prop:properties-of-theta}({\em b})).

Let $J\subseteq\{1,\ldots, n\}$ be the subset so that $C_{\Flag_{12}}\subseteq C_{\Flag_1}$ is given by $u_j=0$ for all
$j\in J$.  For $x\in C_{\Flag_2}$ we then have that

\begin{equation}
\langle \beta_j,x\rangle\rule{0.25cm}{0cm}\mbox{is}\rule{0.25cm}{0cm} 
\left\{
\begin{array}{cl}
\leq 0 & \mbox{if $j\in J$} \\
\geq 0 & \mbox{if $j\not\in J$}
\end{array}
\right.,  
\end{equation} 

and thus, for $j\in J$ and $x\in C_{\Flag_2}$, that the $j$-th coordinate of $(h\circ\expi)(x)$ is $\geq 1$.

For $m\geq 1$ set $p_m:=(\theta\circ h\circ\expi)(x_m)$, and use $(p_m)_i$ to denote the $i$-th coordinate of $p_m$. 
Since $x_m\in C_{\Flag_2}$, we conclude from the discussion above, and the definition of 
$\theta$ (\S\ref{sec:def-of-theta}), that for all $j\in J$, and all $m\geq 1$, $(p_m)_{j}\geq (p_m)_{j+1}$.

By hypothesis the $x_m$ converge to a point $x'\in \overline{C}_{\Flag_1}$.  Let $q\in \Delta_n$ be the 
point corresponding to $x'$ under the homeomorphism $\Delta_n\stackrel{\sim}{\longrightarrow} \overline{C}_{\Flag_1}$
induced by $\psi$.  As above, we use $(q)_i$ to denote the $i$-th coordinate of $q$.  To say that
$x'\in \overline{C}_{\Flag_{12}}$ is the same as saying that $(q)_{j}=(q)_{j+1}$ for all $j\in J$. 

The points $\{p_m\}_{m\geq 1}$ do not necessarily converge to $q$ in $[0,\infty)^n$.  However, by 
the commutativity of the diagram, $\{\psi(p_m)\}_{m\geq 1}$ converges to $\psi(q)$. 

Let $\pi\colon \AA^{m}\longrightarrow \AA^{n}$ be projection
onto the first $n$ coordinates, and $\psi':=\pi\circ\psi$.  We then have that $\{\psi'(p_m)\}_{m\geq 1}$ converges
to $\psi'(q)$ in $\AA^{n}$, and, by our choice of $\alpha_1$, \ldots, $\alpha_n$ at the beginning of the argument,
that the map $\psi'\colon [0,\infty)^n\longrightarrow \AA^{n}$ is a monomial map
given by $\{b_{ij}\}$ satisfying the conditions of Proposition \ref{prop:properties-of-theta}({\em c}). 

Let $i_0$ be the largest index so that the $i_0$-th coordinate of $\psi'(q)$ is zero, or set $i_0=-\infty$ if there is
no such coordinate. From the upper triangular structure of $\psi'$, we conclude (if $i_0\neq -\infty$) that 
$(q)_{i_0}=0$.  Since $q\in \Delta_n$, the inequalities in \eqref{eqn:simplex-boundary} apply, and so $(q)_{j}=0$ for 
all $j\leq i_0$.
Thus, for $j\in J$, $j<i_0$, the boundary condition $(q)_{j}=(q)_{j+1}$ is automatic, and we may concern ourselves
with $j\in J$, $j\geq i_0$.   

By passing to a subsequence, for each $i\in \{1,\ldots, n\}$, we may assume that $\lim_{m\to\infty} (p_m)_i$ exists,
including possibly the limit $\infty$.  We temporarily split $\{1,\ldots, n\}$ into three index sets, 
$\Zsh$, $\Fsh$, and $\Ish$, according to whether or not $\lim_{m\to\infty} (p_m)_i$ is, respectively, zero, 
finite and nonzero, or infinite.   

From the fact that $\lim_{m\to\infty} \psi'(p_m)$ exists and is equal to $\psi'(q)$, and the upper triangular
structure of $\psi'$, we conclude that $i_0=\sup \Zsh$, and that all $j>i_0$ are in $\Fsh$. 
By taking $m$ large enough, or again passing to a subsequence, we may assume that $(p_m)_j>0$ for all $j> i_0$. 
But, again using the upper triangular structure of $\psi'$, since $\psi'(p_m)_j$ is positive and nonzero for
all $j>i_0$, we may solve for each $(p_m)_{j}$ ($j>i_0$) knowing only $\psi'(p_m)$.  Moreover, since this inverse
construction is a continuous map (when all coordinates $>i_0$ are in $\RR_{>0}$), we conclude that
$\lim_{m\to 0} (p_m)_{j} = (q)_{j}$ for all $j\geq i_0$.  (For $j>i_0$, this is the inverse argument above.  For
$j=i_0$, we already know that $\lim_{m\to\infty} (p_m)_{i_0}=0=(q)_{i_0}$.)

For $j\in J$ we have $(p_m)_j\geq (p_m)_{j+1}$ for all $m$.  Thus, for $j\in J$, $j\geq i_0$, by passing to the
limit we conclude that $(q)_{j}\geq (q)_{j+1}$.  The reverse inequality $(q)_{j}\leq (q)_{j+1}$ is part of 
\eqref{eqn:simplex-boundary}.  Thus $(q)_{j}=(q)_{j+1}$ for all $j\in J$, proving ({\em a}). 

To see ({\em b}), given arbitrary flags $\Flag_1$, $\Flag_2$ in $\Delta_X$, set $\Flag_{12}=\Flag_1\cap \Flag_2$.
Choose a maximal flag $\Flag_1'$ containing $\Flag_1$, and let $\Flag'_{12}=\Flag'_{1}\cap \Flag_2$.  
Suppose that $\{x_m\}_{m\geq 1}\in C_{\Flag_2}$ is a sequence of points converging to a point 
$x'\in \overline{C}_{\Flag_1}$.  Since $x'$ is also in $\overline{C}_{\Flag'_1}$, by part ({\em a}) we know 
that $x'\in \overline{C}_{\Flag'_{12}}$. 

Thus $x'\in \overline{C}_{\Flag'_{12}}\cap \overline{C}_{\Flag_1}$.   Although (until we finish proving this theorem)
we do not know how arbitrary closures of cones intersect, we do understand this when the cones being intersected 
are contained in the same maximal cone, since we have an explicit homeomorphism of the closure of each maximal cone
with an $n$-simplex.  In particular, since $\Flag'_{12}\cap \Flag_{1}=\Flag_1'\cap \Flag_2\cap\Flag_1=\Flag_{12}$, 
we know that $\overline{C}_{\Flag'_{12}}\cap \overline{C}_{\Flag_1}=\overline{C}_{\Flag_{12}}$, proving ({\em b}).

For ({\em c}), setting $\Flag_{12}=\Flag_{1}\cap \Flag_{2}$, we already have the inclusion 
$\overline{C}_{\Flag_{12}}\subseteq \overline{C}_{\Flag_1}\cap \overline{C}_{\Flag_2}$.   Any point 
$x'\in \overline{C}_{\Flag_2}$ is a limit of points $\{x_m\}_{m\geq 1}$ with each $x_m\in C_{\Flag_2}$.  Thus, for
any $x'\in \overline{C}_{\Flag_{1}}\cap\overline{C}_{\Flag_2}$, part ({\em b}) 
applies and shows that $x'\in \overline{C}_{\Flag_{12}}$.  \epf

\point At this point we have shown that the $\overline{C}_{\Flag}$, for maximal flags $\Flag$, cover $X_{\geq 0}$,
that each $\overline{C}_{\Flag}$ is an $n$-simplex, and that any two such cones, 
say $\overline{C}_{\Flag_1}$ and $\overline{C}_{\Flag_2}$ fit together, glued along boundary simplices, 
in the same way as the cones $C_{\Flag_1}$ and $C_{\Flag_2}$ (a
statement made precise by Theorem \ref{thm:intersection-of-cones}({\em c}) above). 

Thus one expects that the $\overline{C}_{\Flag}$ glue together to form an $n$-dimensional ball.    

The vector space $N_{\RR}$ has a natural compactification to an $n$-ball $B^{n}$ by adding the ``sphere at $\infty$'' 
(\S\ref{sec:definition-of-BN} below), and perhaps the cleanest way to justify that the $\overline{C}_{\Flag}$ 
glue together to form an $n$-ball, and that the result is homeomorphic to $X_{\geq 0}$ would be to extend a
homeomorphism $N_{\RR}\stackrel{\sim}{\longrightarrow} X_{>0}$ to a 
homeomorphism $B^{n}\stackrel{\sim}{\longrightarrow} X_{\geq 0}$.

The homeomorphism $N_{\RR}\stackrel{\sim}{\longrightarrow} X_{>0}$ provided by $\expi$ does not extend
to a map on the compactification $B^{n}$, see \S\ref{sec:Remarks-on-extension-of-homeomorphism}(\Euler{1}).  
Below we provide a homeomorphism $\Phi\colon N_{\RR}\stackrel{\sim}{\longrightarrow} N_{\RR}$ with the property 
that the composite homeomorphism $\expi\circ\, \Phi$ between $N_{\RR}$ and $X_{>0}$ does 
so extend.
This result is proven in Theorem \ref{thm:top-structure-of-X}.

\bpoint{Definition of $\Phi_{\Flag}$} 
\label{sec:def-of-Phi-F}

For each $k\geq 1$ and $1\leq j\leq k$, define $\phi_{j,k}\colon \RR_{\geq 0}^{k}\longrightarrow \RR_{\geq 0}$
by 

\begin{equation}
\label{eqn:def-of-phi-jk}
\phi_{j,k}(u_1,\ldots, u_k) := \frac{1}{2\pi}\log\left(\frac{1+\sum_{i=1}^{j}u_i}{1+\sum_{i=1}^{j-1}u_i}\right).
\end{equation}

For instance, $\phi_{1,k}(u_1,\ldots, u_k) = \frac{1}{2\pi}\log\left({1+u_1}\right)$ and 
$\phi_{2,k}(u_1,\ldots, u_k)=\frac{1}{2\pi}\log\left(\frac{1+u_1+u_2}{1+u_1}\right)$. 

Let $\Flag=\{\sigma_{i_1}\subset \sigma_{i_2} \subset \cdots \subset \sigma_{i_k}\}$ be a flag in $\Delta_{X}$, 
and $C_{\Flag}\subset N_{\RR}$ the corresponding simplicial cone 
(\S\ref{sec:barycentric-subdivision}).  Each point $x\in C_{\Flag}$ has
a unique expression $x = \sum_{j=1}^{k} u_j \Bsub_{\sigma_j}$ in simplicial coordinates, with each $u_j\geq 0$.
We define $\Phi_{\Flag}\colon C_{\Flag}\longrightarrow C_{\Flag}$ by 

\begin{equation}
\label{eqn:Phi-F-def}
\Phi_{\Flag}(x) = \sum_{j=1}^{k}\phi_{j,k}(u_1,\ldots, u_k) \Bsub_{\sigma_j},
\end{equation}

where, as above, $(u_1,\ldots, u_k)$ are the simplicial coordinates of $x$. 
For instance, if $k=3$, the map is 
$$u_1 \Bsub_{\sigma_1} + u_2\Bsub_{\sigma_2} + u_3\Bsub_{\sigma_3}
\xmapsto{\rule{0.5cm}{0cm}}
\frac{1}{2\pi}\log\left({1+u_1}\right) \Bsub_{\sigma_1} 
+ \frac{1}{2\pi}\log\left(\tfrac{1+u_1+u_2}{1+u_1}\right) \Bsub_{\sigma_2}
+ \frac{1}{2\pi}\log\left(\tfrac{1+u_1+u_2+u_3}{1+u_1+u_2}\right)\Bsub_{\sigma_3}.$$

\bpoint{Proposition/Definition of $\Phi$} 
\label{prop:Phi-F-gluing}

\begin{enumerate}
\item 
Let $\Flag$ be a flag in $\Delta_{X}$ and $C_{\Flag}\subset N_{\RR}$ the corresponding simplicial cone.
Then the map $\Phi_{\Flag}\colon C_{\Flag}\longrightarrow C_{\Flag}$ is a homeomorphism.

\smallskip
\item For each subflag $\Flag'\subset \Flag$, with corresponding face $C_{\Flag'}\subset C_{\Flag}$, 
$$\Phi_{\Flag}|_{C_{\Flag'}} = \Phi_{\Flag'}.$$
\end{enumerate}

If $X$ is proper then as $\Flag$ runs over the flags in $\Delta_{X}$, the cones $C_{\Flag}$ cover $N_{\RR}$. 
By ({\em b}) the maps $\Phi_{\Flag}$ glue together to give a global map 
$\Phi\colon N_{\RR}\longrightarrow N_{\RR}$, which by ({\em a}) is a homeomorphism. 
We note explicitly the defining condition of $\Phi$~: 
for each flag $\Flag$ in $\Delta_{X}$, $\Phi|_{C_{\Flag}} = \Phi_{\Flag}$. 

{\em Proof of the Proposition}. 
In simplicial coordinates the map $\Phi_{\Flag}$ is 

$$
\begin{array}{ccc}
\RR_{\geq 0}^{k} & \xrightarrow{\rule{0.5cm}{0cm}} & \RR_{\geq 0}^{k} \\
(u_1,\ldots, u_k) & \xmapsto{\rule{0.5cm}{0cm}} & 
\left(\phi_{1,k}(u_1,\ldots, u_k),\, \phi_{2,k}(u_1,\ldots, u_k)\ldots, \phi_{k,k}(u_1,\ldots, u_k)\right).
\end{array}
$$

Denoting by $v_1$,\ldots, $v_k$ the coordinates on the target $\RR_{\geq 0}^{k}$, one has 

\begin{equation}
\label{eqn:Phi-inv-def}
u_j = e^{2\pi\sum_{i=1}^{j} v_i} - e^{2\pi\sum_{i=1}^{j-1} v_i}
= \left(e^{2\pi v_j}-1\right) e^{2\pi \sum_{i=1}^{j-1} v_i}.
\end{equation}

The formulas in \eqref{eqn:Phi-inv-def} define a morphism $\RR^{k}_{\geq 0} \longrightarrow \RR^{k}_{\geq 0}$
which is easily checked to be inverse to $\Phi_{\Flag}$, giving ({\em a}).

To see ({\em b}), we suppose that $\Flag = \sigma_{i_1} \subset \sigma_{i_2} \subset \cdots \subset \sigma_{i_k}$, 
and let $\Flag' = \sigma_{i'_1} \subset \cdots \sigma_{i'_{\ell}}$ be a subflag, with $\ell < k$. 

Let $f\colon \left\{1,\ldots, \ell\right\} \longrightarrow \left\{1,\ldots, k\right\}$ 
be the function such that $\sigma_{i'_{j'}} = \sigma_{i_{f(j')}}$ for all $j'\in \{1,\ldots, \ell\}$.
The face $C_{\Flag'}\subset C_{\Flag}$ is defined, in simplicial coordinates on $C_{\Flag}$, by 
the equations $\left\{u_j=0 \st j\not\in \Im(f) \right\}$.  

Let $x$ be a point of $C_{\Flag'}$, and write $x$ in the simplicial coordinates on $C_{\Flag}$~:
$x = (u_1,\ldots, u_k)$.   In the simplicial coordinates on $C_{\Flag'}$, $x$ is then given by
$(u_{f(1)},\ldots, u_{f(\ell)})$.   

\begin{itemize}
\item
If $j\not\in\Im(f)$, then using $u_j=0$ and \eqref{eqn:def-of-phi-jk} we get $\psi_{j,k}(u_1,\ldots, u_k)=\frac{1}{2\pi}\log(1) = 0$.

\smallskip
\item If $j\in \Im(f)$, let $j'$ be the unique element of $\{1,\ldots, \ell\}$ such that $f(j')=j$.  Then, since
$\sum_{i'=1}^{j'} u_{f(i')} = \sum_{i=1}^{j} u_i$, and similarly 
$\sum_{i'=1}^{j'-1} u_{f(i')} = \sum_{i=1}^{j-1} u_i$, we conclude that 
$\psi_{j',\ell}(u_{f(1)},\ldots, u_{f(\ell)}) = \psi_{j,k}(u_1,\ldots, u_k)$. 
\end{itemize}

Together these two statements show that $\Phi_{\Flag}|_{C_{\Flag'}} = \Phi_{\Flag'}$.  \epf

\bpoint{Definition of $B^{n}$}
\label{sec:definition-of-BN}
In this section we recall the standard compactification of a real vector space by the ``sphere at $\infty$''.
Let $N_{\RR}\times N_{\RR}\longrightarrow \RR$, $(x,y)\mapsto x\cdot y$  
be a positive definite symmetric inner product on 
$N_{\RR}$, and let $S_1:=\left\{x\in N_{\RR} \st x\cdot x = 1\right\}$ be the unit sphere.

The map $(0,\infty)\times S_1 \longrightarrow N_{\RR}\setminus\{0\}$ given by $(r,x)\mapsto rx$ is a 
homeomorphism.  We put the standard structure at $\infty$ on $(0,\infty]$, identifying $[0,\infty)$ with $(0,\infty]$ 
via $s\mapsto \frac{1}{s}$ (so that $s=0$ corresponds to $\infty$).  Then $(0,\infty]\times S_1$
can be glued to $N_{\RR}$ along $(0,\infty)\times S_1= N_{\RR}\setminus \{0\}$ to obtain a compact topological space
$B^{n}$, homeomorphic to the $n$-dimensional ball.   We use $S_{\infty}:=\{\infty\}\times S_1\subset B^{n}$ 
for the sphere at $\infty$, and note that there is one point of $S_{\infty}$ for each ray in $N_{\RR}$. 

This construction does not depend on the inner product chosen.  If $S'_1$ is the unit sphere under a different
inner product, radial projection gives a homeomorphism $S_1\stackrel{\sim}{\longrightarrow} S'_1$ inducing a
homeomorphism of the compactifications. 

\bpoint{Theorem} 
\label{thm:top-structure-of-X}
(Topological structure of $X_{\geq 0}$) Let $X$ be a complete toric variety over $\CC$,
$\expi$ and $\Phi$ the maps defined in \S\ref{sec:AMRT} and \S \ref{sec:def-of-Phi-F}
respectively, and put $n=\dim(X)$.

\begin{enumerate}
\item 
The composite homeomorphism $\expi\circ\, \Phi \colon N_{\RR} \stackrel{\sim}{\longrightarrow} X_{>0}$ extends
to a homeomorphism $\overline{\Phi}\colon B_{n} \stackrel{\sim}{\longrightarrow} X_{\geq 0}$.  That is,
we have a commutative diagram 
\PauseEnumerate

\begin{equation}
\xymatrix{
N_{\RR}\ar[rr]^{\Phi}_{\sim}\ar@{^{(}->}[d]\rule[-0.15cm]{0cm}{0.1cm} & & N_{\RR}\ar[rr]^{\expi}_{\sim} & & X_{>0}\ar@{^{(}->}[d]\rule[-0.15cm]{0cm}{0.1cm} \\
B^{n}\ar[rrrr]^{\overline{\Phi}}_{\sim} & & & & X_{\geq 0} \\
}
\end{equation}

\parshape 1 1cm 0.9\textwidth
where the horizontal maps are homeomorphisms and the vertical maps open inclusions.

\ResumeEnumerate
\item $X_{\geq 0}$ is homeomorphic to the $n$-dimensional ball.
\item $X_{\geq 0}$ is a regular cell complex, with a $k$-dimensional cell for each cone $\sigma\in \Delta_{X}$ of
dimension $n-k$. 
\end{enumerate}

\parshape 1 0cm \textwidth

\bpf
({\em a}) Let $\Flag=\{\sigma_1\subset \cdots \subset \sigma_n\}$ be a maximal flag in $\Delta_X$, 
and $C_{\Flag}\subset N_{\RR}$ the corresponding simplicial cone.  We denote by $\tilde{C}_{\Flag}$ the closure
of $C_{\Flag}$ in $B^{n}$.  By radial projection, $\tilde{C}_{\Flag}$ is homeomorphic to the standard 
compactification of $C_{\Flag}$ as an $n$-simplex, which is how we will treat it.
We first show that $\expi\circ\, \Phi_{\Flag}$ extends to a homeomorphism 
$\tilde{C}_{\Flag}\stackrel{\sim}{\longrightarrow} \overline{C}_{\Flag}$. 

Let $U_{\sigma_n}\hookrightarrow \AA^{m}$ be an embedding (given, as usual by choosing a generating set
for $\sigma_n^{\vee}\cap M$), and apply Theorem \ref{thm:C-F-monomial-diagram} to obtain a map $\psi$ making
\eqref{eqn:C-F-closure-factorization} commute.   By 
Theorem \ref{thm:C-F-monomial-diagram}, to show that $\expi\circ\, \Phi_{\Flag}$
extends to a homeomorphism between $\tilde{C}_{\Flag}$ and $\overline{C}_{\Flag}$, it suffices to show
that the composite $\theta\circ h\circ \exp_{\Flag}\circ \Phi_{\Flag}$ extends to a homeomorphism
of $\tilde{C}_{\Flag}$ with $\Delta_n=\Im(\theta)$. 

Combining the definitions of $\Phi_{\Flag}$, $\exp_{\Flag}$, and $\theta$ (see respectively \eqref{eqn:Phi-F-def},
\eqref{eqn:exp-F-formula-2}, and \S\ref{sec:def-of-theta}), we find that, in simplicial coordinates
on $C_{\Flag}$, the map is given by

\begin{equation}
\label{eqn:composite-on-C-F}
\rule{1.2cm}{0cm}
(u_1,\,u_2,\ldots,\, u_n) \mapsto \left(
\frac{1}{1+\sum_{i=1}^{n} u_i}, \, \frac{1+u_1}{1+\sum_{i=1}^{n} u_i},\, \frac{1+u_1+u_2}{1+\sum_{i=1}^{n} u_i},\,
\ldots,\, \frac{1+\sum_{i=1}^{n-1} u_i}{1+\sum_{i=1}^{n} u_i}
\right),
\end{equation}

i.e., the map whose $j$-th coordinate is $\displaystyle\frac{1+\sum_{i=1}^{j-1} u_i}{1+\sum_{i=1}^{n} u_i}$. 

Let us rewrite \eqref{eqn:composite-on-C-F} in terms of the barycentric coordinates on $\tilde{C}_{\Flag}$. 
Suppose that $\xi_0$, \ldots, $\xi_n$ are $\geq 0$, that $\sum_{i=1}^{n} \xi_i=1$, and that $\xi_0\neq 0$.
The corresponding point in $C_{\Flag}$ is given by $u_i=\frac{\xi_i}{\xi_0}$ for $i=1$, \ldots, $n$.
Conversely, given a point in $C_{\Flag}$ with simplicial coordinates $(u_1,\ldots, u_n)$, the
corresponding barycentric coordinates are 
$\xi_0=\frac{1}{1+\sum_{i=1}^{n} u_i}$  and for $j=1$, \ldots, $n$, $\xi_j= \frac{u_j}{1+\sum_{i=1}^{n} u_i}$.

In barycentric coordinates on $\tilde{C}_{\Flag}$, $\xi_0\neq 0$, the map
\eqref{eqn:composite-on-C-F} is 

\begin{equation}
\label{eqn:Phi-bar-simplicial-coords}
(\xi_0,\ldots, \xi_n) \mapsto 
\left(\xi_0,\, \xi_0+\xi_1,\, \xi_0+\xi_1+\xi_2,\, \ldots,\,
\xi_0+\xi_1+\cdots + \xi_{n-1} \rule{0cm}{0.35cm}\right),
\end{equation}

i.e., exactly the parameterization of $\Delta_n$ given in \eqref{eqn:simplex-coordinates}.
This map thus extends (i.e., across $\xi_0=0$) to a homeomorphism 
$\tilde{C}_{\Flag}\stackrel{\sim}{\longrightarrow}\Delta_n$,  and so, using $\psi$, to a homeomorphism 
$$\overline{\Phi}_{\Flag}\colon \tilde{C}_{\Flag}\stackrel{\sim}{\longrightarrow}
\overline{C}_{\Flag}.$$

For a subflag $\Flag'\subset \Flag$, the explicit parameterization \eqref{eqn:Phi-bar-simplicial-coords} shows
that the map above restricted to $\tilde{C}_{\Flag'}$ induces a homeomorphism of $\tilde{C}_{\Flag'}$ with
$\overline{C}_{\Flag'}\subset\overline{C}_{\Flag}$, which we similarly denote by $\overline{\Phi}_{\Flag'}$.

Thus, the extensions on each $\tilde{C}_{\Flag}$ for maximal flags glue together to give a continuous map
$\overline{\Phi}\colon B^{n}\longrightarrow X_{\geq 0}$.   This map is surjective by 
Corollary \ref{cor:structure-of-C-F}({\em d}).  Since both $B^{n}$ and $X_{\geq 0}$ are compact Hausdorff spaces,
to show that $\overline{\Phi}$ is a homeomorphism it then suffices to show that $\overline{\Phi}$ is injective.

Suppose that $x_1$, $x_2\in B^{n}$ are such that $\overline{\Phi}(x_1)=\overline{\Phi}(x_2)$.  
Let $\Flag_1$ and $\Flag_2$ be flags in $\Delta_{X}$ so that $x_1\in \tilde{C}_{\Flag_1}$ and 
$x_2\in \tilde{C}_{\Flag_2}$.  Then $\overline{\Phi}(x_1)\in \overline{C}_{\Flag_1}$ and
$\overline{\Phi}(x_2)\in \overline{C}_{\Flag_2}$.  
Thus, by Theorem \ref{thm:intersection-of-cones}({\em c}), 
$\overline{\Phi}(x_1)= \overline{\Phi}(x_2)$ is contained in 
$\overline{C}_{\Flag_1}\cap \overline{C}_{\Flag_2}=\overline{C}_{\Flag_{12}}$, with $\Flag_{12}=\Flag_1\cap\Flag_2$.

The homeomorphism 
$\overline{\Phi}_{\Flag_1}\colon \tilde{C}_{\Flag_1}\stackrel{\sim}{\longrightarrow}\overline{C}_{\Flag_1}$ 
induces a homeomorphism of the sub-simplex $\tilde{C}_{\Flag_{12}}\subseteq \tilde{C}_{\Flag_1}$ with the sub-simplex
$\overline{C}_{\Flag_{12}}\subseteq \overline{C}_{\Flag_1}$, and thus $x_1\in \tilde{C}_{\Flag_{12}}$.
Similarly $x_2\in \tilde{C}_{\Flag_{12}}$.  Since $\overline{\Phi}_{\Flag_{12}}$ is injective 
(it is a homeomorphism of $\tilde{C}_{\Flag_{12}}$ with $\overline{C}_{\Flag_{12}}$), we conclude that
$x_1=x_2$.  This proves ({\em a}). 

Since $B^{n}$ is an $n$-dimensional ball, part ({\em b}) is an immediate consequence of ({\em a}); 
it was stated separately here for outside reference.

By \cite[Chap. I]{AMRT}, $X_{\geq 0}$ decomposes into $N_{\RR}$-orbits, with an orbit of dimension $k$,
isomorphic to $N_{\RR}/\langle \sigma\rangle$ for each cone $\sigma\in \Delta_{X}$ of dimension $n-k$. 

To show that the cells give $X_{\geq 0}$ the structure of a regular cell complex, the issue is 
to show that the inclusion of each open $k$-cell, thought of as the interior of the $k$-dimensional ball,
extends to a continuous map on the boundary of the ball, which, in addition, is a homeomorphism onto its image.  

Given such a $k$-cell of $X_{\geq 0}$, let $Z$ be the closure of the corresponding complex orbit in $X$ (i.e., 
if the $k$-cell corresponds to $\sigma$, $Z$ is the closure of the complex orbit $O(\sigma)\subset X$).
Then the $k$-cell may be identified with $Z_{>0}$ in $Z_{\geq 0}$.  Since $X$ is proper, $Z$ is proper, and
so applying part ({\em a}) to $Z$ we obtain an extension of $Z_{>0}$ to a homeomorphism of a $k$-dimensional
ball with $Z_{\geq 0}$.  
\epf

\bpoint{Remarks}
\label{sec:Remarks-on-extension-of-homeomorphism}
(\Euler{1}) The inclusion of $N_{\RR}$ into $X_{\geq 0}$ via $\expi$ does not extend to a map from $B^{n}$ 
to $X_{\geq 0}$.  To see why, it is sufficient to pick a maximal flag 
$\Flag = \{\sigma_1\subset \cdots \subset \sigma_n\}$ in $\Delta_{X}$, an embedding 
$U_{\sigma_n}\hookrightarrow \AA^{m}$, and consider the diagram \eqref{eqn:C-F-closure-factorization}
produced by Theorem \ref{thm:C-F-monomial-diagram}.   As in the proof of Theorem \ref{thm:top-structure-of-X}
to understand whether the inclusion of $C_{\Flag}$ in $X_{>0}$ extends to a map 
$\tilde{C}_{\Flag}\longrightarrow X_{\geq 0}$, it suffices to see if the map $\theta\circ h\circ \exp_{\Flag}$
from $C_{\Flag}$ to $[0,1]^n$ has such an extension. 

We note that the question is not whether there some compactification of $C_{\Flag}$ which extends. 
There certainly is, and moreover the extension is an $n$-simplex by 
Theorem \ref{thm:C-F-monomial-diagram} and Corollary \ref{cor:structure-of-C-F}({\em a}).  Rather, the issue is to see 
whether the compactification $\tilde{C}_{\Flag}$ of $C_{\Flag}$ in $B^{n}$ has such an extension.   
As before, we identify $\tilde{C}_{\Flag}$ with the standard compactification of $C_{\Flag}$ as an $n$-simplex.

%
%
%

When $n=2$ the composite map $\theta\circ h\circ \exp_{\Flag}$ is, in simplicial coordinates on $C_{\Flag}$, the
map

\begin{equation}
\label{eqn:extension-counterexample-simplicial}
(u_1,\,u_2) \mapsto \left(e^{-2\pi(u_1+u_2)},\, e^{-2\pi u_2}\right).
\end{equation}

As before let us rewrite \eqref{eqn:extension-counterexample-simplicial} in terms of barycentric coordinates
on $\tilde{C}_{\Flag}$.
For $\xi_0$, $\xi_1$, $\xi_2$ such that $\sum_{i=0}^{2} \xi_i=1$
and each $\xi_i\geq 0$, we have $u_i=\frac{\xi_i}{\xi_0}$, and in these coordinates 
\eqref{eqn:extension-counterexample-simplicial} becomes 

\begin{equation}
\label{eqn:extension-counterexample-barycentric}
(\xi_0,\, \xi_1,\, \xi_2,)\mapsto 
\left(e^{-2\pi\left(\frac{\xi_1+\xi_2}{\xi_0}\right)},\, 
e^{-2\pi\left(\frac{\xi_2}{\xi_0}\right)}
\right).
\end{equation}

The problem is in the second coordinate.  At $(\xi_0,\xi_1,\xi_2)=(0,1,0)$ we obtain the indeterminate form
$e^{-2\pi\left(\frac{0}{0}\right)}$, 
whose value depends on the path in $C_{\Flag}$ approaching $(0,1,0)$. Thus there is no possible continuous
extension of this function to the point $(0,1,0)$ of $\tilde{C}_{\Flag}$. (In contrast, since $\xi_1+\xi_2=1-\xi_0$,
there is no issue in extending the first coordinate function to the boundary;  the limit is the same no matter what the
approaching path.) 

When $n\geq 3$, there are similar issues : only the first coordinate function of $\theta\circ h\circ\exp_{\Flag}$
can be extended to the compactification of $C_{\Flag}$ in $B^{n}$.

(\Euler{2}) The theorem gives a decomposition of $B^{n}$ into $n$-simplices (the $\tilde{C}_{\Flag}$, for $\Flag$
a maximal flag), and for each one, we have an explicit algebraic parameterization of the corresponding $n$-simplex
in $X_{\geq 0}$, as discussed in \S\ref{sec:remark-on-explicit-parameterization}.

(\Euler{3}) Here is a brief discussion of the intuition behind the overall approach to this result.  Fix, as in 
\S\ref{sec:definition-of-BN}, a positive inner product on $N_{\RR}$. For $r>0$ let $S_{r}\subset N_{\RR}$ denote
the sphere of radius $r$, and 
let $S_{\infty}$ denote the sphere at infinity in the compactification $B^{n}$.

For each finite $r>0$, the intersection of $S_{r}$ with the cones in $\Delta_{X}$ gives $S_{r}$ the structure of 
a cell complex.  For instance, if $\sigma\in \Delta_{X}$ is a $k$-dimensional cone, $|\sigma|\cap S_{r}$ is
a $(k-1)$-dimensional cell in the decomposition.  If $X$ is a projective variety, this cell decomposition is
the decomposition into faces of the boundary of a polytope. 

On the other hand, 
the description in \cite[Chap. I]{AMRT} shows that the boundary of $X_{\geq 0}$
also breaks up into cells, although here a $k$-dimensional cone in $\Delta_{X}$ gives an $(n-k)$-dimensional
cone on the boundary, and the order of inclusion between cells and cones is reversed. 
Assuming that $S_{\infty}$ is homeomorphic to the boundary of $X_{\geq 0}$, we thus get a 
cell complex on $S_{\infty}$ dual to the one on $S_{r}$ for finite $r$. 
When $X$ is projective, this is the cell decomposition of the
boundary of the corresponding dual polytope.

Thus, if there is a homeomorphism $B^{n}\stackrel{\sim}{\longrightarrow} X_{\geq 0}$, we would obtain a family
of cell complexes on $S_{r}$ which, on $S_{\infty}$, becomes the dual complex.

The flags of faces of a polytope are in bijection with the flags of faces of its dual.  
This gives a natural bijection between the simplices in the barycentric subdivision of a polytope, and
the simplices in the barycentric subdivision of its dual.

Taking the barycentric subdivision of $N_{\RR}$, it is natural to wonder 
if ({\em i}) the closure in $X_{\geq 0}$ of each $n$-dimensional simplicial cone $C_{\Flag}$ 
has, as boundary, an $(n-1)$-dimensional simplex, 
and ({\em ii}) if these families of $(n-1)$-simplices (on $S_{r}$ and
the boundary $S_{\infty}$) interpolate between the cell complex given by $\Delta_{X}$ and its dual.

The first question is answered by Theorem \ref{thm:C-F-monomial-diagram} and 
Corollary \ref{cor:structure-of-C-F}({\em a}), while Theorem \ref{thm:top-structure-of-X} is, in essence,
an answer to the second question. (In particular, establishing the assumption underlying both questions : 
that there is a homeomorphism $B^{n}\stackrel{\sim}{\longrightarrow} X_{\geq 0}$ which allows us to transfer
the structure of the boundary of $X_{\geq 0}$ to $S_{\infty}$.)

By way of example, let us consider $X=(\PP^1)^3$.  The maximal cones in the fan $\Delta_{X}$ consist of the octants
in $N_{\RR}$, and so the decomposition of $S_{r}$ induced by $\Delta_{X}$ is that of the faces of the octahedron.
On the other hand $X_{\geq 0}$ is a cube, and so the decomposition on $S_{\infty}$, i.e., the boundary of $X_{\geq 0}$,
is that of the faces of a cube.  Figure \ref{fig:barycentric-dual} below shows how the barycentric subdivision
of the octahedron becomes that of the cube ``as $r$ goes to $\infty$''.

\input DualPolytope.tex

\end{document}

%% file: DualPolytope.tex

\newgray{barygray}{0.75}  

\newcommand{\SphereArc}[2]{
\parametricplot{0}{1}{#1 1 t sub \TScale #2 t \TScale \TTrans \MakeUnit \TT }
}

\newcommand{\SphereTriang}[3]{
\pscustom{%
\SphereArc{#1}{#2}
\SphereArc{#2}{#3}
\SphereArc{#3}{#1}
}}

\newcommand{\SphereSmallBaryTri}[3]{
{{
\psset{fillstyle=none,linecolor=barygray}
\SphereArc{#1 #2 \TTrans}{#1 #2 #3 \TTrans \TTrans}
\SphereArc{#1}{#1 #2 #3 \TTrans \TTrans} 
}
\SphereArc{#1}{#1 #2 \TTrans}
\FillSphereTri{#1}{#1 #2 \TTrans}{#1 #2 #3 \TTrans \TTrans}
}}

\newcommand{\SphereBaryTri}[3]{
\SphereSmallBaryTri{#1}{#2}{#3}
\SphereSmallBaryTri{#2}{#1}{#3}
\SphereSmallBaryTri{#2}{#3}{#1}
\SphereSmallBaryTri{#3}{#1}{#2}
\SphereSmallBaryTri{#3}{#2}{#1}
\SphereSmallBaryTri{#1}{#3}{#2}
}

\newcommand{\FillSphereTri}[3]{
{\psset{fillstyle=solid,fillcolor=vlgray,linecolor=vlgray}
\multido{\n=0.04+0.15}{5}{
\SphereTriang{#1 #2 #3 \TTrans \n\space \TScale \TTrans }{#2 #3 #1 \TTrans \n\space \TScale \TTrans }{#3 #1 #2 \TTrans \n\space \TScale \TTrans }
}}
}

\newcommand{\OctBaryTriang}[3]{
\pspolygon(!#1 \TT)(!#1 #2 \TTrans 0.5 \TScale \TT)(!#1 #2 \TTrans #3 \TTrans 0.33333 \TScale \TT)
}

\newcommand{\OctTriang}[3]{
{\psset{linecolor=barygray,fillcolor=vlgray,fillstyle=solid}
\OctBaryTriang{#1}{#2}{#3}
\OctBaryTriang{#2}{#1}{#3}
\OctBaryTriang{#2}{#3}{#1}
\OctBaryTriang{#3}{#2}{#1}
\OctBaryTriang{#3}{#1}{#2}
\OctBaryTriang{#1}{#3}{#2}
}
\pspolygon[linecolor=black,fillstyle=none](!#1 \TT)(!#2 \TT)(!#3 \TT)
}

\newcommand{\CubeBaryTri}[3]{
\pspolygon[fillstyle=solid,fillcolor=vlgray,linecolor=vlgray](!#1 \TT)(!#1 #2 \TTrans 0.5 \TScale \TT)(!#1 #3 \TTrans 0.5 \TScale \TT)
{\psset{fillstyle=none}
\psline[linecolor=barygray](!#1 #2 \TTrans 0.5 \TScale \TT)(!#1 #3 \TTrans 0.5 \TScale \TT)
\psline[linecolor=barygray](!#1 \TT)(!#1 #3 \TTrans 0.5 \TScale \TT)
}
}

\newcommand{\CubeBaryFace}[4]{
\CubeBaryTri{#1}{#2}{#3}
\CubeBaryTri{#2}{#1}{#4}
\CubeBaryTri{#2}{#3}{#4}
\CubeBaryTri{#3}{#2}{#1}
\CubeBaryTri{#3}{#4}{#1}
\CubeBaryTri{#4}{#3}{#2}
\CubeBaryTri{#4}{#1}{#2}
\CubeBaryTri{#1}{#4}{#3}
\pspolygon[linecolor=black,fillstyle=none](!#1 \TT)(!#2 \TT)(!#3 \TT)(!#4 \TT)
}

\RotAngleSetTT{20}{20}

\begin{centering}
\begin{tabular}{ccccccc}
\begin{tabular}{c}
\begin{pspicture}(-1.2,-1.2)(1.2,1.2)  
\OctTriang{-1 0 0}{0 -1 0}{0 0 1}
\OctTriang{1 0 0}{0 1 0}{0 0 1}
\OctTriang{1 0 0}{0 1 0}{0 0 -1}
\OctTriang{-1 0 0}{0 1 0}{0 0 1}
\OctTriang{-1 0 0}{0 1 0}{0 0 -1}
\end{pspicture}
\end{tabular}
& $\rightsquigarrow$ & 
\begin{tabular}{c}
\psset{unit=1.1cm}
\begin{pspicture}(-1.2,-1.2)(1.2,1.2)    
\SphereBaryTri{0 0 -1}{-1 0 0}{0 1 0}  
\SphereBaryTri{1 0 0}{0 1 0}{0 0 1}     
\SphereBaryTri{0 0 1}{-1 0 0}{0 1 0}    
\SphereBaryTri{1 0 0}{0 1 0}{0 0 -1}    
\end{pspicture}
\end{tabular}
& $\rightsquigarrow$ & 
\begin{tabular}{c}
\psset{unit=1.3cm}
\begin{pspicture}(-1.2,-1.2)(1.2,1.2)     
\psset{linecolor=barygray}
\SphereBaryTri{0 0 -1}{-1 0 0 }{0 1 0}  
\SphereBaryTri{1 0 0}{0 1 0}{0 0 1}     
\SphereBaryTri{0 0 1}{-1 0 0}{0 1 0}    
\SphereBaryTri{1 0 0}{0 1 0}{0 0 -1}    
\psset{linecolor=black}
\SphereArc{1 1 1}{-1 1 1}
\SphereArc{-1 1 1}{-1 1 -1}
\SphereArc{-1 1 -1}{1 1 -1}
\SphereArc{1 1 -1}{1 1 1}
\SphereArc{1 1 1}{1 0 1}
\SphereArc{-1 1 1}{-1 0 1}
\SphereArc{1 1 -1}{1 0 -1}
\end{pspicture}
\end{tabular}
& $\rightsquigarrow$ & 
\begin{tabular}{c}
\begin{pspicture}(-1.5,-1.5)(1.5,1.5)          
\CubeBaryFace{1 -1 1}{1 -1 -1}{1 1 -1}{1 1 1}  
\CubeBaryFace{1 1 1}{1 -1 1}{-1 -1 1}{-1 1 1}  
\CubeBaryFace{1 1 1}{-1 1 1}{-1 1 -1}{1 1 -1}  
\end{pspicture}
\end{tabular} \\
\multicolumn{7}{c}{\Fig\label{fig:barycentric-dual}} \\
\multicolumn{7}{c}{\small Example of interpolation between barycentric subdivisions of a polytope and of its dual.} \\
\end{tabular}\\
\end{centering}